\documentclass[12pt]{article}
\usepackage{amsmath}
\usepackage{amssymb}

\usepackage{graphicx}

\def\qed{\hfill$\Box$\par\medskip\par\relax}

\newcommand{\1}[1]{{\mathbb I}{\{#1\}}}

\newcommand{\eps}{\varepsilon}
\newcommand{\Z}{{\mathbb Z}}

\newcommand{\Sph}{{\mathbb S}}

\newcommand{\M}{{\mathcal M}}

\newcommand{\HH}{{\mathcal H}}
\newcommand{\R}{{\mathbb R}}
\newcommand{\RR}{{\mathcal R}_\omega{}}
\newcommand{\CC}{{\mathfrak C}}
\newcommand{\cCC}{{\mathcal C}}

\newcommand{\ZZ}{{\mathcal Z}}
\newcommand{\DD}{{\mathcal D}}
\newcommand{\BB}{{\mathcal B}}
\newcommand{\TT}{{\mathcal T}}

\newcommand{\EEE}{{\mathcal E}}

\newcommand{\e}{\mathbf{e}}

\let\phi=\varphi

\newcommand{\E}{{\mathbf E}}

\newcommand{\gs}{{\mathfrak S}}

\newcommand{\diam}{{\mathop{\rm diam}}}

\newcommand{\IQ}{{\mathbb Q}}

\newcommand{\IP}{{\mathbb P}}

\newcommand{\Po}{{\mathtt P}_\omega}
\newcommand{\Eo}{{\mathtt E}_\omega}

\newcommand{\Eozz}{{\mathtt E}_{\omega,\zeta}}
\newcommand{\Pozz}{{\mathtt P}_{\omega,\zeta}}

\newcommand{\UU}{{\mathcal U}}

\newcommand{\muo}{\mu^\omega}
\newcommand{\nuo}{\nu^\omega}

\newcommand{\nn}{{\mathbf n}_\omega{}}

\newcommand{\exq}[1]{\big\langle #1 \big\rangle_{\!{}_\IQ}}

\newcommand{\exip}[1]{\big\langle #1 \big\rangle_{\!{}_\IP}}

\newcommand{\near}{\stackrel{{\scriptstyle\omega}}{\leftrightarrow}}

\newcommand{\intl}{\int\limits}

\allowdisplaybreaks

\newtheorem{theo}{Theorem}[section]
\newtheorem{lmm}[theo]{Lemma}
\newtheorem{df}[theo]{Definition}
\newtheorem{prop}[theo]{Proposition}
\newtheorem{cor}[theo]{Corollary}
\newtheorem{rem}[theo]{Remark}

\newcommand{\hD}{{\widehat {\mathcal D}}^\omega_H}
\newcommand{\hDl}{{\hat D}_\ell}
\newcommand{\hDr}{{\hat D}_r}
\newcommand{\hF}{{\hat F}^\omega_H}

\newcommand{\tDl}{{\tilde D}_\ell}
\newcommand{\tD}{{\widetilde {\mathcal D}}^\omega_H}

\newcommand{\DDD}{{\mathsf D}}

\title{Knudsen gas in a finite random tube: transport diffusion and
first passage properties}
\author{Francis~Comets$^{1}$ \and
 Serguei~Popov$^{2}$
\and Gunter~M.~Sch\"utz$^{3}$ 
\and Marina Vachkovskaia$^{2}$}

\begin{document}

\maketitle

{\footnotesize
\noindent $^{~1}$Universit{\'e} Paris 7, UFR de Math{\'e}matiques,
case 7012, 2, place Jussieu, F--75251 Paris Cedex 05, France\\
\noindent e-mail: \texttt{comets@math.jussieu.fr},
\noindent url: \texttt{http://www.proba.jussieu.fr/$\sim$comets}

\smallskip
\noindent $^{~2}$Department of Statistics,
Institute of Mathematics, Statistics and Scientific Computation,
University of Campinas--UNICAMP,
rua S\'ergio Buarque de Holanda 651, CEP 13083--859,
Campinas SP, Brazil\\
\noindent e-mails: \texttt{popov@ime.unicamp.br},
\texttt{marinav@ime.unicamp.br}\\
\noindent urls:
\texttt{http://www.ime.unicamp.br/$\sim$popov},
\texttt{http://www.ime.unicamp.br/$\sim$marinav}

\smallskip
\noindent $^{~3}$Forschungszentrum J\"ulich GmbH,
Institut f\"ur Festk\"orperforschung,
D--52425 J\"ulich, Deutschland\\
\noindent e-mail: \texttt{G.Schuetz@fz-juelich.de}, \\
\noindent url:
\texttt{http://www.fz-juelich.de/iff/staff/Schuetz\_G/}

}

\begin{abstract}
We consider transport diffusion in a stochastic billiard in a random
tube which is elongated in the direction
of the first coordinate (the tube axis). Inside the random tube,
which is stationary and ergodic, non-interacting particles
move straight with constant speed. Upon hitting the tube walls, they
are reflected randomly, according to the cosine law: the density of
the outgoing direction is proportional to the cosine of the angle
between this direction and the normal vector. Steady state transport
is studied by introducing an open tube
segment as follows: We cut out a large finite segment of the tube
with segment boundaries perpendicular
to the tube axis. Particles which leave this piece through the
segment boundaries disappear from the system. 
Through stationary injection of particles at one boundary of the
segment a steady state with non-vanishing 
stationary particle current
is maintained. We prove (i) that in the thermodynamic limit of an
infinite open piece the coarse-grained density profile inside the
segment is linear, and (ii) that the transport diffusion
coefficient obtained from the ratio of stationary current and
effective boundary density gradient equals the diffusion coefficient
of a tagged particle in an infinite tube. Thus we prove Fick's law
and equality of transport diffusion and self-diffusion 
coefficients for quite generic rough (random) tubes. 
We also study some properties of the crossing time and compute the
Milne extrapolation length in dependence on the shape
of the random tube.
 \\[.3cm]\textbf{Keywords:} cosine law, Knudsen random walk,
random medium, self-diffusion coefficient, transport diffusion
coefficient, random walk in random environment
\\[.3cm]\textbf{AMS 2000 subject classifications:}
 60K37. Secondary: 37D50, 60J25
\end{abstract}

\section{Introduction}
\label{s_intro}

Diffusion in stationary states may be encountered either in
equilibrium, where no macroscopic mass or energy fluxes are present
in a system of many diffusing particles, or away from equilibrium,
where diffusion is often driven by a density gradient between two
open segments of the surface that encloses the space in which
particles diffuse. In equilibrium states, one is interested in the
\emph{self-diffusion} coefficient $\DDD_{\text{self}}$, as given
by the mean-square displacement (MSD) of a tagged particle. This
quantity, also called tracer diffusion coefficient, can be
measured using e.g.\ neutron scattering, NMR or direct video 
imaging in the case of colloidal particles. In gradient-driven
non-equilibrium steady states, there is a particle flux between the
boundaries which is proportional to the density gradient. This factor
of proportionality is the so-called transport or collective diffusion
coefficient $\DDD_{\text{trans}}$.

Often these two diffusion coefficients cannot be measured
simultaneously under concrete experimental conditions and the
question arises whether one can infer
knowledge about the other diffusion coefficient, given one of them.
Generally, in dense systems these diffusion coefficients depend in a
complicated fashion on the interaction between the diffusing
particles. In the case of diffusion in microporous media, e.g.\ in
zeolites, however, the mean free path of the particles is of the
order of the pore diameter or even larger. Then diffusion is
dominated by the interaction of particles with the
pore walls rather than by direct interaction between particles. In
this dilute so-called Knudsen regime 
neither $\DDD_{\text{self}}$ nor $\DDD_{\text{trans}}$
depend on the particle density anymore, but are just given by the
low-density limits of these two quantities. One then expects
$\DDD_{\text{self}}$ and $\DDD_{\text{trans}}$ 
to be equal. This assumption is a fundamental input into the
interpretation of many experimental data, see e.g.~\cite{HK} for an
overview of diffusion in condensed matter systems.

Not long ago this basic tenet has been challenged by Monte-Carlo
simulation of Knudsen diffusion in pores with fractal pore walls
\cite{MC1,MC2,MC3}. The authors of these (and further) studies
concluded that self-diffusion depends on the surface
roughness of a pore, while transport diffusion is independent of it.
In other words, the authors of  \cite{MC1,MC2,MC3} argue that even in
the low density limit, where the gas particle are independent of each
other and interact only with the pore walls,
$\DDD_{\text{self}} \neq \DDD_{\text{trans}}$, with a dependence of
$\DDD_{\text{self}}$ on the details of the pore walls that 
$\DDD_{\text{trans}}$ does not exhibit.
This counterintuitive numerical finding was quickly questioned on
physical grounds and contradicted by 
further simulations~\cite{RZBK} which give approximate equality of
the two diffusion coefficients. These controversial results
gave rise to a prolonged debate which finally led to the consensus
that indeed both diffusion coefficients should agree for the Knudsen
case~\cite{ZRBCK}. It has remained open though
whether these diffusion coefficients are generally exactly equal or
only approximately to a degree depending on the details of the
specific setting.

A physical argument put forward in~\cite{ZRVCDBK} suggests general
equality. To see this one imagines the following
\textit{gedankenexperiment}. 
Imagine one colours in a equilibrium setting of
many non-interacting particles
some of these particles without changing their properties. At some
distance from this colouring region the colour is removed. Then these
coloured particles experience a density gradient just as ``normal''
particles in an open system with the same pore walls would. Since the
walls are essentially the same and the properties
of coloured and uncoloured particles are the same, the statistical
properties of the ensemble of trajectories remain unchanged. Hence
one expects any pore roughness to have the same effect on diffusion,
irrespective of whether
one consider transport diffusion or self-diffusion. Notice, however,
that this microscopic argument, while intuitively appealing, is far
from rigorous. First, the precise conditions under which the
independence of the diffusion
coefficients on the pore surface is supposed to be valid, is not
specified. This is more than a technical issue since one may easily
construct surface properties 
leading to non-diffusive behaviour (cf.~\cite{CPSV2,MVW}).
Second, there is no obvious microscopic interpretation or unique
microscopic definition of the transport diffusion coefficient for
arbitrary surface structures. $\DDD_{\text{trans}}$ is a genuinely
macroscopic quantity and a proof of equality between
$\DDD_{\text{trans}}$ and $\DDD_{\text{self}}$
(which is naturally microscopically defined through the asymptotic
long-time behaviour of the MSD) requires some further work and new
ideas. One needs
to establish that on large scales the Knudsen process converges to
Brownian motion (which then also gives $\DDD_{\text{self}}$).
Moreover, in order to compare $\DDD_{\text{trans}}$ and
$\DDD_{\text{self}}$ one needs a precise macroscopic definition of
$\DDD_{\text{trans}}$ which is independent of microscopic properties
of the system.

The first part of this programme is carried out in~\cite{CPSV2}. 
There we proved the quenched 
invariance principle for the horizontal
projection of the particle's position using the 
method of considering
the environment viewed from the particle.
This method is useful in a number of models related to
Markov processes in a random environment, 
cf.\ e.g.~\cite{DFGW,FSS,M}.
The aim of this paper is to solve the second problem of defining
$\DDD_{\text{trans}}$
and proving equality with $\DDD_{\text{self}}$. As in~\cite{CPSV2}
we consider a random tube to model pore roughness. 
In contrast to~\cite{CPSV2}, we now have to consider
tubes of finite extension along the tube contour and introduce open
segments at the ends of the tube. Doing this rigorously then
clarifies some of the salient assumptions underlying the equality of
$\DDD_{\text{trans}}$ and $\DDD_{\text{self}}$. Naturally, since we
are in the dilute gas limit, there is no dependence on the particle
density in either of the two diffusion constants. This obvious point
has not been controversial and will not be stressed below.

We note that we define $\DDD_{\text{trans}}$ through
stationary transport in an open system since this is accessible
experimentally as well as numerically in Monte Carlo simulation.
Indeed, in the literature that gave rise to the
controversy that we address here, this way of defining
$\DDD_{\text{trans}}$ is used, albeit in a non-rigorous fashion.
Sticking to this experimentally motivated
  setting we shall give below a precise definition
that can be used to prove rigorously that under rather generic
circumstances $\DDD_{\text{trans}}=\DDD_{\text{self}}$, which means
that both diffusion constants depend on the pore surface in the same
way. As pointed out above, this equality 
 is expected from independence of the particles and the invariance
principle for the process and its time-reversed.
However, 
 we could not find a general result applying here, and moreover,
as it turns out, the proof is not entirely
trivial. 
There are some technical difficulties to overcome because
the quenched invariance principle of
Definition~\ref{def_invariance_principle} below is not very
``strong'' (there is no uniformity assumption on the speed of
convergence as a function of the initial conditions)
 and the jumps of
the embedded discrete-time billiard are not uniformly bounded. Let
us mention here that it is generally difficult to obtain stronger
results in the above sense, since the corrector technique, generally
used in the proof of quenched central limit theorems for reversible
Markov processes in random environment, is still not sufficiently
well understood.

To further illuminate the contents of our results we point out that
in a bulk system the equality of the self-diffusion coefficient and
the transport diffusion coefficient for the spread of
\textit{equilibrium} density fluctuations in an infinite system may
be taken for granted in 
the case of particles that have no mutual interaction. Hence another
way of stating the main conclusion of our work is the assertion that
the transport diffusion coefficient as defined here in a stationary
\textit{far-from-equilibrium} setting coincides with the usual
equilibrium transport diffusion coefficient. 

We also address finite-size effects coming from the fact that
we are dealing with diffusion in a finite, open geometry. 
This causes deviations from bulk results for first-passage-time
properties if a tagged particle starts its motion close to one
boundary. In particular, we compute the permeation time and the
Milne extrapolation length that characterizes the survival time of a
particle injected at a boundary.

As a final introductory remark, it is worth noting that the case of
Knudsen gas with the cosine reflection law (which is the model
considered in this paper) is particularly easy to analyse
because the stationary state can be 
written in an explicit form, cf.\
Theorem~\ref{t_stat_measure}. As explained below, this is
related to the following facts: (i) there is no interaction between
particles, (ii) for random billiard (i.e., a motion of only one
particle in a closed domain) with the cosine reflection law the
stationary measure is quite explicit, as shown in~\cite{CPSV1}.
Similar questions are much more complicated
when the explicit form of the stationary state
is not known. This is the general situation for non-equilibrium
steady states. We refer to e.g.\ the model of~\cite{BeOl} 
(a chain of
coupled oscillators) where one resorts to a bound on the entropy
production.

This paper is organized in the following way. In
Section~\ref{s_defin} we define the infinite random tube, and then
introduce the process we call random billiard. In
Section~\ref{s_transport}, we then consider a gas of independent
particles with absorption/injection in a finite piece of the random
tube, and we formulate our results on the stationary measure for
that gas and on the transport diffusion coefficient. In
Section~\ref{s_perm}, we go on to formulate first passage time
 results that concern
exit from and crossing of the finite tube by a tagged particle.
The remaining part of the paper is devoted to the
proof of our results. In Section~\ref{s_pr_prelim} we mainly use
the reversibility of the process to obtain several technical facts
used later. In Section~\ref{s_pr_steady} we prove the result
on the stationary measure of the Knudsen gas in the finite tube.
 Section~\ref{s_pr_trans} contains the proofs of the results
related to the transport diffusion coefficient, and in
Section~\ref{s_pr_perm} we prove the results related to the crossing
of the finite tube.

\section{General notations and main results}
\label{s_main_results}

Naively the transport diffusion coefficient in tube
direction~$x$ may be defined through the diffusion equation for the 
probability density $\partial_t P(x,t) = \partial_x ( D(x) 
\partial_x P(x,t) )$, where a possible $x$-dependence may originate
from a spatial inhomogeneity of the tube.
Denote by~$J$ the
particle current in the system; assuming stationarity with
a probability density $P^*(x)$ one has 
$J  = D(x) \partial_x P^*(x)$.
With fixed external densities~$P^+$ at $x=L$ 
and~$P^-$ at $x=0$ one finds by integration 
$J = \DDD_{\text{trans}}\vartheta$
with density gradient $\vartheta = (P^+ - P^-)/L$
and $\DDD_{\text{trans}}^{-1} = 1/L \int_0^L dx D^{-1}(x)$. 
By measuring the current and the boundary densities
one can thus obtain the transport diffusion coefficient without
having to determine the local quantity $D(x)$.
This result,
however, implies knowledge of the local coarse-grained boundary 
densities~$P^\pm$ to be able to make any comparison with
$\DDD_{\text{self}}$. In a 
real experimental setting as well as for a given microscopic
model these boundary densities~$P^\pm$ are difficult
to obtain. In particular, there is no well-defined prescription where
precisely on a microscopic scale these boundary
quantities should be measured.
We circumvent the problem of computing these quantities
from microscopic considerations
by considering the total number of particles in the tube rather
than local properties of the boundary region of the tube. Together
with proving a large-scale linear density profile in a stationary
open random tube, one may then infer the macroscopic density
gradient, see the definition~(\ref{defdtrans}) below. Thus one
obtains a macroscopic definition of
the transport diffusion coefficient which is independent of
microscopic details of the model.

\subsection{Definitions of the random tube and the random billiard}
\label{s_defin}
In order to fix ideas in a mathematically rigorous form we 
first recall some notations from~\cite{CPSV2}.

Let us formally define the random tube in~$\R^d$, $d\geq 2$.
In this paper, $\R^{d-1}$ will always stand for the linear
subspace of $\R^d$ which is perpendicular to the first coordinate
vector~$\e$,
we use the notation~$\|\cdot\|$ for the Euclidean norm in~$\R^d$
or~$\R^{d-1}$. For $k\in \{d-1,d\}$
let $\BB(x,\eps)=\{y\in\R^k:\|x-y\|<\eps\}$ be the open
$\eps$-neighborhood of $x\in\R^k$. Define
 $\Sph^{d-1}=\{y\in\R^d:\|y\|=1\}$ to be the unit sphere in~$\R^d$.
Let
\[
\Sph_h = \{w\in\Sph^{d-1}: h\cdot w > 0\}
\]
be the half-sphere looking in the direction~$h$.
 For $x\in\R^d$, sometimes it will be convenient to write
$x=(\alpha,u)$,
being~$\alpha$ the first coordinate of~$x$ and $u\in\R^{d-1}$; then,
$\alpha=x\cdot\e$,
and we write $u=\UU x$, being~$\UU$ the projector on~$\R^{d-1}$. 
Fix some positive
constant~${\widehat M}$, and define
\begin{equation}
\label{def_Lambda}
 \Xi = \{u\in \R^{d-1} : \|u\| \leq {\widehat M}\}.
\end{equation}

Let~$A$ be an open connected domain in~$\R^{d-1}$ or~$\R^d$.
We denote by~$\partial A$ the boundary of~$A$ and 
by~$\bar A = A\cup \partial A$ the closure of~$A$.

The random tube is viewed as 
a stationary and ergodic process $\omega=(\omega_\alpha,
\alpha\in\R)$,
where~$\omega_\alpha$ is a subset of~$\Xi$; cf.~\cite{CPSV2}
for a more detailed definition. 
We denote by~$\IP$ the law of this process; sometimes 
we will use the shorthand notation $\exip{\cdot}$ for the
expectation with respect to~$\IP$. 
With a slight abuse of notation, we denote also by
\[
 \omega = \{(\alpha,u)\in\R^d : u\in \omega_\alpha\}
\]
the random tube itself, where the billiard lives. 
Intuitively, $\omega_\alpha$ is the ``slice'' obtained 
by crossing~$\omega$ with the hyperplane 
$\{\alpha\}\times \R^{d-1}$. 
We will assume that the domain $\omega$ is defined in such a way
that it is an open subset of $\R^d$, and that it is connected.
We write also~$\bar\omega$ for the closure of~$\omega$.
In order to define the random billiard correctly,
following~\cite{CPSV1}, 
throughout this paper we suppose that $\IP$-almost
surely~$\partial\omega$ is a
$(d-1)$-dimensional surface satisfying the Lipschitz condition.
This means that for any~$x\in\partial\omega$ there exist~$\eps_x>0$,
an affine isometry ${\mathfrak I}_x : \R^d\to\R^d$, a function
$f_x:\R^{d-1}\to\R$
such that
\begin{itemize}
\item $f_x$ satisfies Lipschitz condition, i.e., there exists a
constant~$L_x>0$
such that $|f_x(z)-f_x(z')| < L_x\|z-z'\|$ for all $z,z'$;
\item ${\mathfrak I}_x x = 0$, $f_x(0)=0$, and
\[
 {\mathfrak I}_x(\omega\cap\BB(x,\eps_x)) = \{z\in\BB(0,\eps_x) : 
                           z^{(d)} > f_x(z^{(1)},\ldots,z^{(d-1)})\}.
\]
\end{itemize}
Roughly speaking, Lipschitz condition implies that any boundary point
can be ``touched'' by a piece of a cone which lies fully inside the
tube. This in its turn ensures that the (discrete-time) process
cannot remain in a small neighborhood of some boundary point for very
long time; in Section~2.2 of~\cite{CPSV1} one can find an example
of a non-Lipschitz domain where the random billiard behaves in an
unusual way.  

We keep the usual notation $dx, dv, dh, \ldots$ for the
$(d-1)$-dimensional Lebesgue measure
on~$\Xi$ (usually restricted to $\omega_\alpha$ for
some~$\alpha$) or Haar measure
on~$\Sph^{d-1}$. We write $|A|$ for the $k$-dimensional Lebesgue
measure in case $A\subset \R^k$,
and Haar measure in case $A\subset \Sph^{d-1}$.
Also, we denote by~$\nuo$ the $(d-1)$-dimensional Hausdorff measure
on $\partial\omega$; since the boundary is Lipschitz, one obtains
that~$\nuo$  is locally finite (cf.\ the proof of Lemma~3.1
in~\cite{CPSV1}).

We assume additionally that the boundary of $\IP$-a.e.~$\omega$ is
$\nuo$-a.e.\ \emph{continuously} differentiable, and we 
denote by $\RR\subset\partial\omega$ the set of boundary points
where~$\partial\omega$ is continuously differentiable.

To avoid complications when cutting a (large) finite piece of the
infinite random tube,
we assume that there exists a constant $\widetilde M$ such that 
for $\IP$-almost all environments~$\omega$ we have the following: 
for any $x,y\in\omega$ with $|(x-y)\cdot\e|\leq 1$
there exists a path connecting $x,y$ that lies fully inside~$\omega$ 
and has length at most~$\widetilde M$.

For all $x\in\RR$, let us define the normal
vector $\nn(x)\in\Sph^{d-1}$ pointing inside the domain~$\omega$.

We say that $y\in\bar\omega$ is \emph{seen from} $x\in\bar\omega$ if
there exists $h\in\Sph^{d-1}$ and~$t_0>0$ such that $x+th\in\omega$
for all $t\in (0,t_0)$
and $x+t_0 h = y$. 
Clearly, if~$y$ is seen from~$x$ then~$x$ is seen from~$y$,
and we write ``$x \near y$'' when this occurs. 

Next, we construct the Knudsen random walk
(KRW) $(\xi_n, n=0,1,2,\ldots)$, which is a discrete time Markov
process on~$\partial\omega$,
cf. Section~2.2 of~\cite{CPSV1}.
It is defined through its transition density~$K$:
for $x,y\in\partial\omega$
\begin{equation}
\label{def_K}
 K(x,y) = \gamma_d
\frac{\big((y-x)\cdot\nn(x)\big)\big((x-y)\cdot\nn(y)\big)}{\|x-y\|^{
d+1}}
\1{x,y\in\RR, x \near y},
\end{equation}
where $\gamma_d = \big(\int_{\Sph_{\e}} h\cdot\e\, dh\big)^{-1}$
is the normalizing constant, and~$\1{\cdot}$ stands for the indicator
function.
This means that, being $\Po,\Eo$ the quenched (i.e., with
fixed~$\omega$)
probability and expectation, for any $x\in\RR$ and any measurable
$B\subset \partial\omega$
we have
\[
 \Po[\xi_{n+1}\in B \mid \xi_n=x] = \intl_B K(x,y)\, d\nuo(y).
\]
We also refer to the Knudsen random walk as the random walk with
cosine reflection law, since it is elementary to obtain
from~\eqref{def_K} that the density
of the outgoing direction is proportional to the cosine of the angle
between this direction and the normal vector.
\begin{rem}
\label{rem_not_to_infty}
 In fact, in the general setting of~\cite{CPSV1}, for unbounded
domains, one has to consider
the following possibility: at some moment the particle 
chooses the outgoing direction in such a way that, moving in this
direction, it never hits the 
boundary of the domain again, thus going directly to the infinity.
However,
it is straightforward to see that, since
$\omega\subset\R\times\Xi$, 
in our situation $\Po$-a.s.\ this cannot happen.
\end{rem}

 It is immediate to obtain from~\eqref{def_K} that
$K(\cdot,\cdot)$ is symmetric (that is, $K(x,y)=K(y,x)$ for all
$x,y\in\partial\omega$);
for both the discrete- and continuous-time processes
this leads to some nice reversibility properties, exploited
in~\cite{CPSV1, CPSV2}. 
Clearly, $K$ depends on~$\omega$ as well, but we usually do not
indicate this in the notations in order to keep them simple. 
Also, let us denote by~$K^n(\cdot,\cdot)$ the $n$-step transition
density; clearly, one obtains that~$K^n$ is symmetric too for
any~$n\geq 1$. 

Now, we define the Knudsen stochastic billiard (KSB)
$\big((X_t,V_t),t\geq 0\big)$, which is the main object of study in
this paper. 
First, we do that for the process starting on the
boundary~$\partial\omega$
from the point~$x_0\in \partial\omega$. Let
$x_0=\xi_0,\xi_1,\xi_2,\xi_3,\ldots$
be the trajectory of the random walk, and define
\[
 \tau_n = \sum_{k=1}^n \|\xi_k-\xi_{k-1}\|.
\]
Then, for $t\in[\tau_n,\tau_{n+1})$, define 
\[
 X_t=\xi_n+(\xi_{n+1}-\xi_n)\frac{t-\tau_n}{\|\xi_{n+1}-\xi_n\|}.
\]
In Proposition~2.1 of~\cite{CPSV1} it was shown that, provided that
the boundary satisfies the Lipschitz condition, we have
$\tau_n\to\infty$ $\Po$-a.s., and so~$X_t$
is well-defined for all $t\geq 0$.
The quantity~$X_t$ stands for the position of the particle at
time~$t$; since it is not a Markov process by itself, we define also
the c\`adl\`ag version of the motion direction at time~$t$:
\[
 V_t = \lim_{\eps \downarrow 0}\frac{X_{t+\eps}-X_t}{\eps},
\]
observe that $V_t\in \Sph^{d-1}$. 
Recall also another notation from~\cite{CPSV1}:
for $x\in\omega$, $v\in\Sph^{d-1}$, define (with the convention
$\inf\emptyset=\infty$)
\[
{\mathsf h}_x(v) = x+v\inf\{t>0 : x+tv \in \partial\omega\} \in
\partial\omega\cup\{\infty\},
\]
so that ${\mathsf h}_x(v)$ is the next point where the particle hits 
the boundary when starting at the location~$x$ with the
direction~$v$. Of course, we can define also the stochastic billiard
starting from the interior of~$\omega$ by specifying its initial
position~$x_0$ and initial direction~$v_0$: the particle starts at
the position~$x_0$ and moves
in the direction~$v_0$ with unit speed until hitting the boundary
at the point ${\mathsf h}_{x_0}(v_0)$; then, the previous 
construction is applied, being ${\mathsf h}_{x_0}(v_0)$ the starting
boundary point. We denote by $\Po^{x,v}$ the (quenched) law of KSB in
the tube~$\omega$ starting from~$x$ with the initial direction~$v$.

Consider the rescaled projected trajectory ${\hat Z}_t^{(s)} =
s^{-1/2}X_{st}\cdot\e$ of KSB.
\begin{df}
\label{def_invariance_principle}
We say that the quenched invariance principle holds for the Knudsen
stochastic billiard 
in the infinite random tube if there exists a positive
constant~${\hat \sigma}$ such that, 
for $\IP$-almost all~$\omega$, 
for any initial conditions $(x_0,v_0)$ such that ${\mathsf
h}_{x_0}(v_0)\in\RR$,
the rescaled trajectory ${\hat \sigma}^{-1}{\hat
Z}^{(s)}_\cdot(\omega)$
weakly converges to the Brownian motion as $s\to\infty$.
\end{df}

Also, for some of our results we will have to make more assumptions
on the geometry of the random tube. Consider the following

\medskip
\noindent
\textbf{Condition T.}
\begin{itemize}
\item[(i)] There exists a positive constant~$\bar\eps$ and a
continuous function $\bar\phi : \R\mapsto \R^d$ such that
\[
 \inf_{\substack{t\in\R\\ x\in\R^d\setminus\omega}}
 \|\bar\phi(t)-x\| \geq \bar\eps,
 \quad \lim_{t\to-\infty}\bar\phi(t)\cdot\e =-\infty, 
 \quad \lim_{t\to\infty}\bar\phi(t)\cdot\e =\infty.
\]
\item[(ii)] In the case $d\geq 3$, we assume that there
exist~$N,r_1>0$ such that for all $x,y\in\RR$ 
with $|(x-y)\cdot\e|\leq 2$ there exists $n\leq N$ such that
$K^n(x,y)\geq r_1$.
\item[(iii)]
In the case $d=2$, we assume that
\[
 \sup\{|(x-y)\cdot\e| : x,y\in\RR, x\near y\} < \infty 
           \qquad \text{$\IP$-a.s.}
\]
\end{itemize}

\begin{rem}
\label{rem_almost_all}
 From the fact that $\omega\subset\R\times\Xi$ and $\nuo$-almost
all points of~$\partial\omega$
belong to~$\RR$, it is straightforward to obtain that for
Lebesgue$\times$Haar-almost all
$(x,v)\in\omega\times\Sph^{d-1}$ we have ${\mathsf h}_x(v)\in\RR$
(see Lemma~3.2~(i) of~\cite{CPSV1}).
\end{rem}

\begin{rem}
\label{rem_CLT}
 In the paper~\cite{CPSV2} we prove that, if the second moment 
of the projected jump length with respect to
the stationary measure for the environment seen from the particle 
is finite (which is true for $d\geq 3$, but not always for $d=2$), 
then under certain additional
conditions (related to Condition~T of the present paper), 
the quenched invariance principle holds for the Knudsen
stochastic billiard in the infinite random tube, cf.\ Theorem~2.2,
Propositions~2.1 and~2.2 of~\cite{CPSV2}. 
Let us comment more on the above Condition~T:
\begin{itemize}
\item In~\cite{CPSV2}, instead of the ``uniform D\"oblin condition''
(ii), we assumed a more explicit (although a bit more technical)
Condition~P, which implies that (ii) holds (see Lemma~3.6
of~\cite{CPSV2}). In fact, in the proof of the quenched invariance
principle the technical condition of~\cite{CPSV2} is used only
through the fact that it implies the uniform D\"oblin condition.
 \item The assumption we made for $d=2$ may seem to be too
restrictive. However, is it only a bit more restrictive that the
assumption that the random tube does not contain an  infinite
straight cylinder. As it was shown in Proposition~2.2
of~\cite{CPSV2}, if 
the random tube contains an infinite straight cylinder, then
the averaged second moment of the projected jump length is infinite
in dimension~$2$, and so the (quenched) invariance principle cannot
be valid.
\end{itemize}
\end{rem}

\subsection{Gas of independent particles and evaluation of the
transport diffusion coefficient}
\label{s_transport}

\begin{figure}
 \centering
\includegraphics{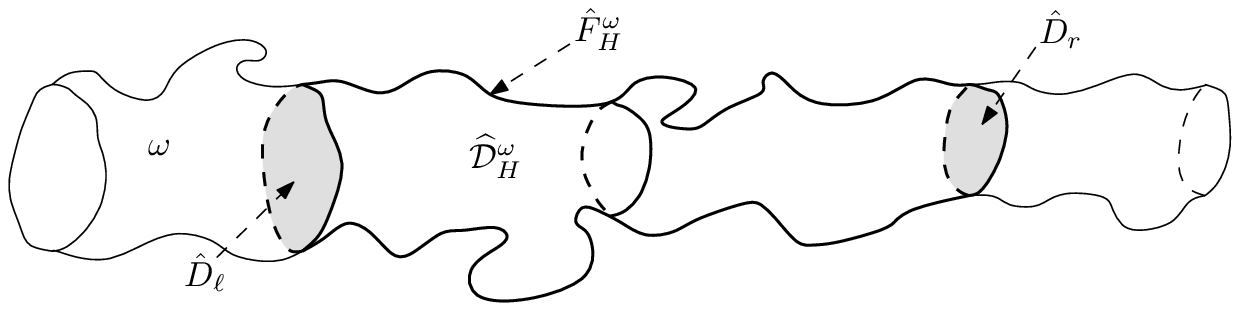}
\caption{On the definition of finite tube~$\hD$}
\label{f_piece_of_tube}
\end{figure}

Now, let us introduce the notations specific to this paper. 
Consider a positive number~$H$ (which is typically
supposed to be large); denote by~$\hD$ the part 
of the random tube~$\omega$ which lies between~$0$ and~$H$:
\[
 \hD = \{z\in\omega : z\cdot\e \in [0,H]\}.
\]
Denote also 
\begin{align*}
 \hF &= \{x\in\partial\omega : x\cdot\e \in (0,H)\},\\
 \hDl &= \{0\}\times \omega_0,\\
 \hDr &= \{H\}\times \omega_H,
\end{align*}
so that $\partial\hD = \hF \cup \hDl \cup \hDr$ (see
Figure~\ref{f_piece_of_tube}). 
Observe that~$\hD$ can, in fact, consist of several separate pieces,
namely, one big piece between~$0$ and~$H$, 
and possibly several small
pieces near the left and the right ends
(we suppose that $H\geq\widetilde M$, so that there could not be two
or more big pieces). It can be easily seen that those small pieces
have no influence on the definition of the transport diffusion
coefficient; for notational convention, we still allow~$\hD$ 
to be as described above.

Then, we consider a gas of independent particles in~$\hD$, described
as follows. There is 
usual reflection on~$\hF$; any particle which hits $\hDl\cup \hDr$,
disappears.
In addition, for a given $\lambda>0$,
new particles are injected in~$\hDl$ with
intensity~$(\gamma_d|\Sph^{d-1}|)^{-1}\lambda$
per unit surface area.
Every newly injected particle chooses the initial 
direction at random according to the cosine law.
In other words, the injection in~$\hDl$ is given by an independent
Poisson process in $\hDl \times \Sph_{\e}$
with intensity $|\Sph^{d-1}|^{-1} \lambda | \e \cdot u|\, dx\, du$. 
\begin{rem}
\label{rem_cosine_chosen}
The choice of the cosine law for the injection of new particles is
justified
by Theorem~2.9 of~\cite{CPSV1}: for the KSB in a finite domain, the
long-run empirical
law of intersection with a $(d-1)$-dimensional manifold is cosine.
One may think of the following
situation: the random tube is connected from its left side~$\hDl$ to
a very large reservoir containing the Knudsen gas in the stationary
regime; then, the particles cross~$\hDl$
with approximately cosine law (at least on the time scale when the
density of the particles in the big reservoir remains unaffected by
the outflow through the tube).
In Section~\ref{s_pr_steady}
(proof of Theorem~\ref{t_stat_measure}) we use this kind
of argument to obtain a rigorous characterization of the steady
 state of this gas.
\end{rem} 

We now consider this gas in the stationary regime. 
Let $\Xi_{[a,b]}:=[a,b]\times\Xi$, and let~$\M(a,b)$ be the
mean number of
particles in $\hD\cap\Xi_{[a,b]}$, in a fixed
environment~$\omega$.

In Theorem~\ref{t_dens_grad} below we shall see that there
exists a constant~$\vartheta$ such that
\[
 \lim_{m\to\infty}\limsup_{H\to\infty}
 \max_{j=1,\ldots,m}\Bigg|
 \frac{\M\big(\frac{(m-j)H}{m},\frac{(m-j+1)H}{m}
\big)}
     {H/m} - \frac{\vartheta(j-1/2)}{m}\Bigg| = 0,
\]
which means that, after coarse-graining, 
the particle density profile
is asymptotically linear. The above quantity~$\vartheta$ is called
the (rescaled) density gradient.

We define also the current $J^\omega_H$ as the mean number of
particles absorbed in~$\hDr$
per unit of time, and let the rescaled current be defined as
\[
 J=\lim_{H\to\infty}H J^\omega_H.
\]
Then, consistently with the discussion in the beginning of this
section, 
the \emph{transport diffusion coefficient} $\DDD_{\text{trans}}$ is
defined by
\begin{equation}
\label{defdtrans}
 \DDD_{\text{trans}} = \frac{J}{\vartheta}.
\end{equation}
Now, suppose that the quenched invariance principle with
constant~${\hat\sigma}$ holds for the stochastic billiard.  
Our goal is to prove that $\DDD_{\text{trans}}$ is equal to the
\emph{self-diffusion coefficient}
$\DDD_{\text{self}}:={\hat\sigma}^2/2$. 
To this end, we prove the following two results. 
First, we prove that the coarse-grained density profile
is indeed linear:
\begin{theo}
\label{t_dens_grad}
Suppose that the quenched invariance principle holds. Then,
 for any $\eps'>0$ there exists~$m$ such that $\IP$-a.s.
\begin{equation}
\label{eq_dens_grad}
 \limsup_{H\to\infty}
\max_{j=1,\ldots,m}\Bigg|
\frac{\M\big(\frac{(m-j)H}{m},\frac{(m-j+1)H}{m}
\big)}
     {H/m} - \frac{\lambda (j-1/2)}{m}\exip{|\omega_0|}\Bigg| 
        < \eps'
\end{equation}
\end{theo}
Then, we calculate the limiting current:
\begin{theo}
\label{t_current}
Suppose that the quenched invariance principle holds with
constant~${\hat\sigma}$, and assume also that Condition~T
holds. Then, we have $\IP$-a.s.
\begin{equation}
\label{eq_current}
\lim_{H\to\infty}H J^\omega_H = \frac{1}{2}\lambda{\hat\sigma}^2
\exip{|\omega_0|}.
\end{equation}
\end{theo}

Some remarks are in place that illustrate the significance
of the above theorems.
 Theorem~\ref{t_dens_grad} means that $\vartheta = \lambda
\exip{|\omega_0|}$,
and using also Theorem~\ref{t_current}, we obtain that
$\DDD_{\text{trans}}=\DDD_{\text{self}}$.
At the same time it becomes clear that such a statement 
can be true only asymptotically since in a finite
open tube one has to expect finite size corrections of the mean
particle number. These corrections
may, in fact, depend strongly on the microscopic shape of the tube
near the open
boundaries. This implies that in experiments on real spatially
inhomogeneous systems some care 
has to be taken as to what is measured as macroscopic density
gradient.
Notice that with Theorem~\ref{t_current} we also prove 
Fick's law for diffusive
transport of matter in the random Knudsen stochastic billiard. 
Since the velocity of the particles does not change at
collisions with the tube walls, mass transport is proportional
to energy transport. In this interpretation Theorem~\ref{t_current}
implies Fourier's law for heat conduction, see e.g.~\cite{BeOl, GG}
for recent work on other processes.

For a function $g\in\cCC[0,\infty)$ and $a\in\R$, denote
\begin{equation} 
\label{eq:f000999}
 \wp_a(g) = \inf\{t\geq 0 : g(t)-g(0)=a\}.
\end{equation}
As mentioned in the introduction, in the proof of
Theorems~\ref{t_dens_grad} and~\ref{t_current} we use the explicit
form of the steady state for the Knudsen gas in the random tube with
injection from one side. 
Let us formulate the following theorem:
\begin{theo}
\label{t_stat_measure}
{~}
\begin{itemize}
\item[(i)] For the Knudsen gas with
absorption/injection in $\hDr\cup\hDl$ (as before, with intensity
$(\gamma_d|\Sph^{d-1}|)^{-1}\lambda$ per unit surface area)
the unique stationary state is Poisson point process 
in $\hD\times \Sph^{d-1}$ with intensity~$\lambda|\Sph^{d-1}|^{-1}$.
\item[(ii)]
For the gas with injection in~$\hDl$ only, the unique stationary
distribution of the particle configuration is given by a Poisson
point process in~$\hD\times\Sph^{d-1}$ with intensity measure 
\[
 \lambda |\Sph^{d-1}|^{-1} \Po^{(\alpha,u),-h}[\wp_{-\alpha}(X\cdot
e)<\wp_{H-\alpha}(X\cdot\e)]
  \,d\alpha\, du\, dh.
\]
\end{itemize}
Also, in both cases, for any initial configuration the process
converges to the stationary state described above.
\end{theo}
Of course, the above result is not quite unexpected. 
It is well known that independent systems have Poisson invariant
distributions (with the single particle invariant measure for 
Poisson
intensity), let us mention e.g.~\cite{D} (Section~VIII.5)
and~\cite{L}. Still, we decided to include the proof of this theorem
because (as far as we know), it does not directly follow from any of
the existing results available in the literature.

\subsection{Crossing time properties}
\label{s_perm}
\begin{figure}
 \centering
\includegraphics{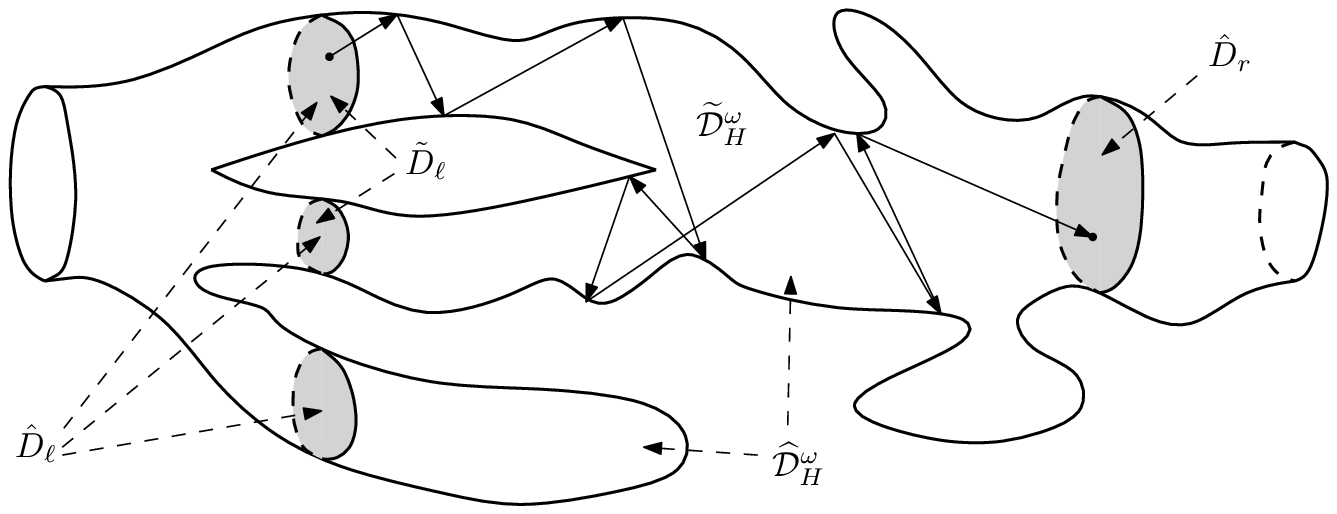}
\caption{On the definition of $\tDl,\tD$, and the event~$\CC_H$
 (a trajectory crossing the tube is shown)}
\label{f_permeation}
\end{figure}

Let us introduce some more notations for the finite random tube.
We denote by ${\tilde\omega}_0$ the set of points of~$\omega_0$,
 from where the particle can reach~$\hDr$ by a path which stays
within~$\hD$ and set $\tDl:=\{0\}\times{\tilde\omega}_0$ (see
Figure~\ref{f_permeation}), and
let~$\tD\subset\hD$ be the corresponding finite tube. Since we
are going to study now how long a tagged particle stays inside the
tube and how it crosses (i.e., goes to the right boundary without
going back to the left boundary), the idea is to
inject it in a place from where it can actually do it. 
Our interest is then in certain first-passage properties, 
in particular, the total life time of the particle inside~$\tD$
(i.e., the time until the particle first exits~$\tD$) and the
permeation time which the particle needs to first exit~$\tD$ at the
end of the tube segment ``opposite'' to that
where it was injected, i.e., after crossing the tube.

So, suppose that one particle is injected (uniformly) at random
at~$\tDl$ into the tube~$\tD$ (that is, the starting location has
the uniform distribution in~$\tDl$, and the direction is chosen
according to the cosine law), and let us denote by~$\CC_H$ the event
that it crosses the tube without going back to~$\tDl$, i.e.,
$\CC_H=\{\tau(\hDr)<\tau^+(\tDl)\}$ (here, $\tau$ and~$\tau^+$ are,
respectively, entrance and hitting times for the discrete-time
process, see~\eqref{def_entrance} and~\eqref{def_hitting} for the
precise definitions). 
  Also, define~$\TT_H$ 
to be the total lifetime of the particle, i.e., if~$X_t$ is the
location of the particle at time~$t$, then 
$\TT_H=\min\{t>0:X_t\in \tDl\cup\hDr\}$.

First, we calculate the asymptotic behaviour of the quenched and
annealed (averaged) expectation of~$\TT_H$:
\begin{theo}
\label{t_exp_perm}
 Suppose that the quenched invariance principle holds
with constant~${\hat\sigma}$. We have
\begin{align}
\lim_{H\to\infty} \frac{1}{H}\Eo\TT_H &= 
    \frac{\gamma_d|\Sph^{d-1}|\exip{|\omega_0|}}
          {2|{\tilde\omega}_0|}  \qquad \text{$\IP$-a.s.,}
\label{q_exp_perm}\\
\lim_{H\to\infty} \frac{1}{H} \exip{\Eo\TT_H} &= 
  \frac{1}{2}\gamma_d|\Sph^{d-1}|\exip{|\omega_0|}
     \exip{|{\tilde\omega}_0|^{-1}}.
\label{a_exp_perm}
\end{align}
\end{theo}
Observe that Condition~T~(i) implies that $|{\tilde\omega}_0|$ is
bounded away from~$0$, and so
$\exip{|{\tilde\omega}_0|^{-1}}<\infty$. 
At this point we  remind the reader that 
here and in the next theorem
the expected ``times'' are actually expected lengths of flight,
related through the corresponding times through the trivial generic
relation
\textit{length}$\,=\,$\textit{velocity}$\times$\textit{time}.
In our Knudsen gas we always assume unit velocity $v=1$ so that 
times can be identified with the appropriate lengths.

To elucidate the physical significance of Theorem \ref{t_exp_perm} 
we observe that for usual Brownian motion the expected lifetime
$T(z_0)$ of particle in an interval $[0,L]$ is given by 
$T(z_0) = z_0 (L-z_0)/(2D)$, where~$z_0$ is the starting position
and~$D$ is the diffusion coefficient. So, in particular, for a
particle starting at the boundary $z_0 = 0$ (or at $z_0=L$) the
expected life time is~$0$. However, in a microscopic model of
diffusion in a finite open system, 
 this result cannot be expected to be generally valid
because of a positive probability that a particle which starts at
$z_0=0$ would escape through the other boundary at~$L$. Often it is
found empirically that the expected life time
can be approximated by 
\begin{equation}
\label{Milne}
 T(\tilde{z}_0) = \frac{\tilde{z}_0 (\tilde{L}-\tilde{z}_0)}{2D}
\end{equation}
with an effective shifted coordinate $\tilde{z}_0 = z_0 + \lambda_M$
and effective interval length $\tilde{L} = L + 2\lambda_M$. The
empirical shift length $\lambda_M$ is known as Milne extrapolation
length~\cite{Case67}, for a recent application to 
diffusion in carbon nanotubes see~\cite{HDMC}.
 From the definition~(\ref{Milne}) one can see that the life time of
a particle starting at the origin $z_0=0$ allows for the computation
of the Milne extrapolation length through the asymptotic relation
\[
 \lim_{L\to\infty} \frac{T(\lambda_M)}{L} = \frac{\lambda_M}{2D}
\]
provided the diffusion coefficient~$D$ is known.

In a physical system the Milne extrapolation length depends on 
molecular details of the gas such as type of molecule or 
temperature, but in a Knudsen gas also on the tube surface.
In our model the properties of the gas are encoded in the unit
velocity~$v=1$ of the particles. Observe now that the quantity
$T(\lambda_M)$ corresponds to~$\Eo\TT_H$ in our setting.
Hence, by identifying $H=L$ and using $D=\hat{\sigma}^2/2$,
 Theorem~\ref{t_exp_perm} 
furnishes us with the dependence of the Milne extrapolation length 
on the tube properties through
\begin{align}
\lambda_M &= 
    \frac{\gamma_d|\Sph^{d-1}|\exip{|\omega_0|}}
          {2|{\tilde\omega}_0|}  \hat{\sigma}^2 \qquad \text{$\IP$-a.s.,}
\label{Milne2a}\\
 \exip{\lambda_M} &= 
  \frac{1}{2}\gamma_d|\Sph^{d-1}|\exip{|\omega_0|}
     \exip{|{\tilde\omega}_0|^{-1}} \hat{\sigma}^2.
\label{Milne2b}
\end{align}
Interestingly, $\lambda_M$ depends only on very few generic 
properties of the random tube.

The next result relies on Theorem~\ref{t_current}, so we need to
assume a stronger condition on the geometry of the tube.
\begin{theo}
\label{t_perm}
Let us suppose that the quenched invariance principle is valid
with~${\hat\sigma}$, and assume that Condition~T holds. 
For the asymptotics of the probability of crossing, we have
\begin{align}
 \lim_{H\to\infty} H\Po[\CC_H] =
\frac{\gamma_d|\Sph^{d-1}|{\hat\sigma}^2\exip{|\omega_0|}}
          {2|{\tilde\omega}_0|}
\qquad \text{$\IP$-a.s.,}
\label{q_perm_prob}\\
\lim_{H\to\infty} H\exip{\Po[\CC_H]} = \frac{1}{2}
\gamma_d|\Sph^{d-1}|{\hat\sigma}^2\exip{|\omega_0|}
  \exip{|{\tilde\omega}_0|^{-1}}.
\label{a_perm_prob}
\end{align}
For the quenched behaviour of the conditional expectations, we have,
$\IP$-a.s.
\begin{align}
 \lim_{H\to\infty} \frac{1}{H^2}\Eo(\TT_H\mid \CC_H)
         &= \frac{1}{3{\hat\sigma}^2},
\label{q_cond_crossing}\\
 \lim_{H\to\infty} \frac{1}{H}\Eo(\TT_H\1{\CC_H})
   & = \frac{\gamma_d|\Sph^{d-1}|\exip{|\omega_0|}}
          {6|{\tilde\omega}_0|},
\label{q_crossing}\\
\lim_{H\to\infty} \frac{1}{H}\Eo(\TT_H\1{\CC_H^c})
   & = \frac{\gamma_d|\Sph^{d-1}|\exip{|\omega_0|}}
          {3|{\tilde\omega}_0|},
\label{q_not_crossing}
\end{align}
and for the annealed ones
\begin{align}
 \lim_{H\to\infty} \frac{1}{H^2}\exip{\Eo(\TT_H\mid \CC_H)}
         &= \frac{1}{3{\hat\sigma}^2},
\label{a_cond_crossing}\\
 \lim_{H\to\infty} \frac{1}{H}\Eo(\TT_H\1{\CC_H})
   & = \frac{1}{6}\gamma_d|\Sph^{d-1}|\exip{|\omega_0|}
          \exip{|{\tilde\omega}_0|^{-1}},
\label{a_crossing}\\
\lim_{H\to\infty} \frac{1}{H}\Eo(\TT_H\1{\CC_H^c})
   & = \frac{1}{3}\gamma_d|\Sph^{d-1}|\exip{|\omega_0|}
          \exip{|{\tilde\omega}_0|^{-1}}.
\label{a_not_crossing}
\end{align}
\end{theo}

As one sees from Theorems~\ref{t_exp_perm} and~\ref{t_perm}, all our
annealed results in fact say that one can interchange the limit as
$H\to\infty$ with integration with respect to~$\IP$. We still
decided to include these results (even though they are technically
not difficult) because, in models related to random environment, it
is frequent that the annealed behaviour 
differs substantially from the quenched behaviour.

One may find it interesting to observe that, by~\eqref{q_crossing}
and~\eqref{q_not_crossing}
\[
 \frac{\Eo(\TT_H\1{\CC_H^c})}{\Eo(\TT_H\1{\CC_H})}\to 2
 \quad \text{ as }H\to\infty.
\]
To obtain another interesting consequence of our results,
let us suppose now that $\IP$-a.s.\ the random tube is such that 
we have $|{\tilde\omega}_0|=|\omega_0|$.
Observe that, by Jensen's inequality, it holds that 
\[
\exip{|\omega_0|}\exip{|\omega_0|^{-1}}\geq 1
\]
(and the inequality is strict if the distribution of~$|\omega_0|$
is nondegenerate), so ``roughness'' of the tube makes the quantities
$\exip{\Eo\TT_H}$ and $\exip{\Eo(\TT_H\1{\CC_H})}$ 
increase. In other
words, these quantities as well as the Milne correlation length
are minimized on the tubes with constant
section (which, by the way, do not have to be necessarily ``straight
cylinders''!).

The remaining part of the paper is devoted to the proofs of our
results, and,
as mentioned in the introduction, it is organized in the
following way.
In Section~\ref{s_pr_prelim} we obtain several auxiliary
results related to hitting of sets by the random billiard.
In Section~\ref{s_pr_steady} we obtain the explicit form of
the stationary measure of the Knudsen gas in the finite tube~$\hD$ 
by using the corresponding result from~\cite{CPSV1} about the
stationary distribution of one particle in a finite domain.
Then, in Section~\ref{s_pr_trans}, we apply the results of
Sections~\ref{s_pr_prelim} and~\ref{s_pr_steady} to 
obtain the explicit form of
 the transport diffusion coefficient.
Finally, in Section~\ref{s_pr_perm} we use Little's theorem to
prove the results related to the crossing time of the random tube.

\section{Some preliminary facts: hitting times and estimates on
the crossing probabilities}
\label{s_pr_prelim}
We need first to prove several auxiliary facts for random billiard 
in arbitrary finite domains. As in~\cite{CPSV1}, let~$\DD$ be a
bounded domain with Lipschitz and a.e.\ continuously differentiable
boundary. We keep the notation~$\Po$ to denote the law of our
processes, and we still use~$\nuo$ to denote the $(d-1)$-dimensional
Hausdorff measure on the boundary~$\partial\DD$.
Consider a Markov chain~$\bar\xi$ on~$\partial\DD$, which
has a transition density~$\bar K$ with the property
${\bar K}(x,y)={\bar K}(y,x)$ for all $x,y\in\partial\DD$.
Observe that the Knudsen random walk~$\xi$ has the above property,
but we need to formulate the next results in a slightly more general
framework, since we shall need to apply them to some other processes
built upon~$\xi$.
Let us introduce the notations
\begin{align}
 \tau(B) &= \min\{n\geq 0 : \bar\xi_n\in B\}, \label{def_entrance}\\
 \tau^+(B) &= \min\{n\geq 1 : \bar\xi_n\in B\} \label{def_hitting}
\end{align}
for the entrance and the hitting time of~$B\subset\partial\DD$.
Also, for measurable~$B\subset\partial\omega$ such that 
$0<\nuo(B)<\infty$ we shall write
\[
 \Po^B[\cdot] = \frac{1}{\nuo(B)}\intl_B\Po^x[\cdot]\,d\nuo(x),
\]
so that~$\Po^B$ is the law for the process starting from the
uniform distribution on~$B$.

Taking advantage of the reversibility of the process~$\xi$,
we prove the following
\begin{lmm}
\label{l_hitting_revers}
 Consider two arbitrary measurable sets $B,F\subset\partial\DD$
such that $B\cap F = \emptyset$.
\begin{itemize}
\item[(i)] Suppose that $\nuo(B),\nuo(F)\in (0,+\infty)$. For any
$F'\subset F$, we have
\begin{align}
 \lefteqn{\Po^B[\xi_{\tau(F)}\in F'\mid \tau(F)<\tau^+(B)]}
\nonumber\\
&  =\frac{1}{\nuo(B)\Po^B[\tau(F)<\tau^+(B)]}
\intl_{F'}\Po^y[\tau(B)<\tau^+(F)]\, d\nuo(y)
\nonumber\\
 &  =\frac{1}{\nuo(F)\Po^F[\tau(B)<\tau^+(F)]}
\intl_{F'}\Po^y[\tau(B)<\tau^+(F)]\, d\nuo(y).
\label{hitting_rev1}
\end{align}
\item[(ii)] Suppose that $\nuo(B)\in (0,+\infty)$. For any
$B',B''\subset B$, we have
\begin{align}
\lefteqn{\intl_{B'} \Po^x[\xi_{\tau^+(B)}\in B'', \tau^+(B)<\tau(F)]
\,d\nuo(x)} \nonumber\\
 & =  \intl_{B''} \Po^x[\xi_{\tau^+(B)}\in B', \tau^+(B)<\tau(F)]
\,d\nuo(x).
\label{hitting_rev2}
\end{align}
\end{itemize}
\end{lmm}
One immediately obtains the following consequence of
Lemma~\ref{l_hitting_revers}~(ii):
\begin{cor}
\label{c_trans_dens}
 For any $B,F\subset\partial\DD$ such that $B\cap F =
\emptyset$ and $\nuo(B)\in (0,+\infty)$, we
have the following.
\begin{itemize}
\item[(i)] For $x,y\in B$, let us define the conditional (on the
event $\{\tau^+(B)<\tau(F)\}$) transition density ${\bar
K}_{B,F}(x,y)$:
\[
 \Po^x[\xi_{\tau^+(B)}\in B''\mid \tau^+(B)<\tau(F)]
 = \intl_{B''} {\bar K}_{B,F}(x,y)\,d\nuo(y).
\]
Then, we have 
\[
 \Po^x[\tau^+(B)<\tau(F)] {\bar K}_{B,F}(x,y) = 
   \Po^y[\tau^+(B)<\tau(F)] {\bar K}_{B,F}(y,x),
\]
that is, the random walk conditioned to return to~$B$ without
hitting~$F$ is reversible with the reversible measure~$\nuo_{B,F}$
defined by
\[
 \frac{d\nuo_{B,F}}{d\nuo}(x) = \Po^x[\tau^+(B)<\tau(F)].
\]
\item[(ii)] In particular (take $F=\emptyset$ in the previous part) 
the random
walk observed at the moments of successive visits to~$B$ is
reversible with the reversible measure~$\nuo$.
\end{itemize}
\end{cor}

\noindent
\textit{Proof of Lemma~\ref{l_hitting_revers}.}
Abbreviate for the moment $U:=\DD\setminus(B\cup F)$.
First, write using the fact that~${\bar K}$ is symmetric
\begin{align*}
 \Po^B[\tau(F)<\tau^+(B)] 
&= \sum_{n=1}^\infty \Po^B[\tau(F)=n,\tau^+(B)>n]\\
 &= \sum_{n=1}^\infty \intl_B \frac{d\nuo(x_0)}{\nuo(B)}
 \intl_{U^{n-1}} d\nuo(x_1)
\ldots d\nuo(x_{n-1}) \\
 &\qquad\qquad\times \intl_F d\nuo(x_n)
 {\bar K}(x_0,x_1)\ldots {\bar K}(x_{n-1},x_n)\\
 &= \frac{\nuo(F)}{\nuo(B)}\sum_{n=1}^\infty \intl_F
\frac{d\nuo(x_n)}{\nuo(F)}
 \intl_{U^{n-1}} d\nuo(x_{n-1})
\ldots d\nuo(x_1) \\
 &\qquad\qquad\times \intl_B d\nuo(x_0)
 {\bar K}(x_n,x_{n-1})\ldots {\bar K}(x_1,x_0)\\
&= \frac{\nuo(F)}{\nuo(B)} \sum_{n=1}^\infty
\Po^F[\tau(B)=n,\tau^+(F)>n]\\
&= \frac{\nuo(F)}{\nuo(B)} \Po^F[\tau(B)<\tau^+(F)].
\end{align*}
Then, similarly
\begin{align*}
 \lefteqn{\Po^B[\xi_{\tau(F)}\in F'\mid \tau(F)<\tau^+(B)]}
\\
 &=\frac{1}{\Po^B[\tau(F)<\tau^+(B)]} 
  \sum_{n=1}^\infty \Po^B[\tau(F)=n,\xi_{\tau(F)}\in F',
\tau^+(B)>n]\\
 &=\frac{\nuo(B)}{\nuo(F)\Po^F[\tau(B)<\tau^+(F)]}
 \sum_{n=1}^\infty \intl_B \frac{d\nuo(x_0)}{\nuo(B)}
 \intl_{U^{n-1}} d\nuo(x_1)
\ldots d\nuo(x_{n-1}) \\
 &\qquad\qquad\qquad\qquad\qquad\qquad\qquad\times \intl_{F'}
d\nuo(x_n)
 {\bar K}(x_0,x_1)\ldots {\bar K}(x_{n-1},x_n)\\
&= \frac{1}{\nuo(F)\Po^F[\tau(B)<\tau^+(F)]} 
  \intl_{F'}\Po^y[\tau(B)<\tau^+(F)]\, d\nuo(y),
\end{align*}
so~\eqref{hitting_rev1} is proved.

Let us prove~\eqref{hitting_rev2}. Analogously to the previous
computation, we write
\begin{align*}
 \lefteqn{ \intl_{B'} \Po^x[\xi_{\tau^+(B)}\in B'',
\tau^+(B)<\tau(F)]
\,d\nuo(x)}\\
 &= \intl_{B'}d\nuo(x)
\sum_{n=1}^\infty \Po^x[\xi_{\tau^+(B)}\in B'',
\tau^+(B)=n,\tau(F)>n]\\
 &= \sum_{n=1}^\infty \intl_{B'} d\nuo(x_0)
 \intl_{U^{n-1}} d\nuo(x_1)
\ldots d\nuo(x_{n-1}) \\
 &\qquad\qquad\times \intl_{B''}
d\nuo(x_n)
 {\bar K}(x_0,x_1)\ldots {\bar K}(x_{n-1},x_n)\\
&= \sum_{n=1}^\infty \intl_{B''} d\nuo(x_n)
 \intl_{U^{n-1}} d\nuo(x_{n-1})
\ldots d\nuo(x_1) \\
 &\qquad\qquad\times \intl_{B'}
d\nuo(x_0)
 {\bar K}(x_n,x_{n-1})\ldots {\bar K}(x_1,x_0)\\
&= \intl_{B''} \Po^x[\xi_{\tau^+(B)}\in B', \tau^+(B)<\tau(F)]
\,d\nuo(x),
\end{align*}
and~\eqref{hitting_rev2} is proved. This concludes the proof of
Lemma~\ref{l_hitting_revers}. 
\qed

\medskip

Next, we recall the Dirichlet's principle:
\begin{prop}
\label{p_Dirichlet}
Consider $B,F\subset \partial\DD$ with $B\cap F=\emptyset$
and $\nuo(B)\in (0,+\infty)$, 
and denote
 ${\hat h}(x) = \Po^x[\tau(F)<\tau(B)]$ (so that, in particular, 
${\hat h}(x)=0$
for all~$x\in B$ and ${\hat h}(x)=1$ for all~$x\in F$). Define
\[
\HH = \{h : h(x)\in[0,1], h(x)=0 \text{ for all } x\in B, 
h(x)=1 \text{ for all } x\in F\}.
\]
Then
\begin{equation}
\label{variational}
 2\nuo(B)\Po^B[\tau(F)<\tau^+(B)] = \EEE({\hat h},{\hat h})
   = \min_{h\in \HH}\EEE(h,h),
\end{equation}
where
\begin{equation}
\label{def_Dirichlet_form}
\EEE(h,h) = \intl_{(\partial\DD)^2}{\bar K}(x,y)(h(x)-h(y))^2 
       \, d\nuo(x)\, d\nuo(y).
\end{equation}
\end{prop}

\noindent
\textit{Proof.}
 For the proof, we refer to the discrete case, e.g.~Proposition 3.8
in~\cite{AF}, and observe that the proof applies to the
space-continuous case, using that, on 
general spaces, harmonicity in the analytic sense and in the
probabilistic sense are equivalent notions by~\cite{C}. Indeed,
minimizers~$h$ of the Dirichlet form are  harmonic in the analytic
sense, i.e., there are in the kernel 
of the form (see~(2.10) in~\cite{C}), 
though the left-hand side of~(24) is the value of $\EEE(h,h)$
when~$h$ is harmonic in the probabilistic 
sense, i.e., the expectation of the process at some exit time (see
Theorem~2.7 in~\cite{C}) with the appropriate boundary conditions.
\qed

\medskip

Now, we go back to the Knudsen random walk in the random
tube~$\omega$. Recall that~$K^n$ stands for the the $n$-step
transition density of KRW, and that we have $K^n(x,y)=K^n(y,x)$ for
all~$x,y$.

Let us define for an arbitrary $A\subset\R$
\[
{\tilde F}^\omega(A) = \{x\in\partial\omega : x\cdot
e \in A\}.
\]
In case~$A$ is an interval, say, $A=[a,b)$, we write
${\tilde F}^\omega[a,b)$ instead of ${\tilde F}^\omega([a,b))$.
There is the following apriory bound on the size of the jump of the
random billiard: there exists
a constant~${\tilde\gamma}_1>0$, depending only
on~${\widehat M} = \diam(\Xi)/2$
and the dimension, such that for $\IP$-almost all~$\omega$
\begin{equation}
\label{bound_d>=4}
\Po[|(\xi_1-\xi_0)\cdot\e|\geq u \mid\xi_0=x] \leq {\tilde\gamma}_1
u^{-(d-1)},
\end{equation}
for all $x\in\partial\omega$, $u\geq 1$, see formula~(54)
of~\cite{CPSV2}. Moreover, using~\eqref{bound_d>=4}, for any $n\geq
1$ it is straightforward to obtain that, for some
${\tilde\gamma}^{(n)}_1>0$
\begin{equation}
\label{bound_n_d>=4}
\Po[|(\xi_n-\xi_0)\cdot\e| \geq u \mid \xi_0=x] \leq 
{\tilde\gamma}^{(n)}_1 u^{-(d-1)},
\end{equation}
for all $x\in\partial\omega$, $u\geq 1$ (also, without restriction
of generality, we can assume that ${\tilde\gamma}^{(n)}_1$ is
nondecreasing in~$n$).

Now, with the help of the above formula we prove the following
result:
\begin{lmm}
\label{l_separated}
 For any $n\geq 1$ there exists ${\tilde\gamma}^{(n)}_2>0$ such that
for all $u\geq 1$ and $a\in\R$ we have
\begin{equation}
\label{eq_separated}
 \intl_{{\tilde F}^\omega(-\infty,a)}d\nuo(x)
  \intl_{{\tilde F}^\omega(a+u,\infty)}d\nuo(y)
   K^n(x,y) \leq {\tilde\gamma}^{(n)}_2 u^{-(d-1)}.
\end{equation}
\end{lmm}

\noindent
\textit{Proof.}
Abbreviate $V=\{a+u\}\times\omega_{a+u}$. The main idea is the
following: if at some step the Knudsen random walk jumped from
some point of ${\tilde F}^\omega(-\infty,a+u)$ to ${\tilde
F}^\omega[a+u,\infty)$, it must cross~$V$, so the probability of
such a jump is the same as the probability of the jump to~$V$ in the
semi-infinite tube with the boundary 
${\tilde F}^\omega(-\infty,a+u)\cup V$. So, we obtain
\begin{align*}
\lefteqn{\intl_{{\tilde F}^\omega(-\infty,a]}d\nuo(x)
  \intl_{{\tilde F}^\omega[a+u,\infty)}d\nuo(y)
   K^n(x,y)}\qquad\\
 &= \intl_{{\tilde F}^\omega(-\infty,a]}
  \Po^x[\xi_n\in {\tilde F}^\omega[a+u,\infty)]\,d\nuo(x)\\
 &\leq \intl_{{\tilde F}^\omega(-\infty,a]}
 \Po^x\Big[\bigcup_{k=1}^n\{\xi_k\cdot\e \geq a+u, 
   \xi_j\cdot\e < a+u \text{ for all }j<k\}\Big]\, d\nuo(x)\\
&\leq \intl_{{\tilde F}^\omega(-\infty,a]}d\nuo(x_0)
 \sum_{k=1}^n \intl_{({\tilde F}^\omega(-\infty,a+u))^{k-1}}
    d\nuo(x_0)\ldots d\nuo(x_{k-1})\\
&\qquad\qquad\qquad\qquad\qquad\qquad
 \intl_{{\tilde F}^\omega[a+u,\infty)}d\nuo(x_k)
 K(x_0,x_1)\ldots K(x_{k-1},x_k)\\
 &\leq \intl_{{\tilde F}^\omega(-\infty,a]}d\nuo(x)
 \intl_V d\nuo(y)
 \big(K(x,y)+K^2(x,y)+\cdots+K^n(x,y)\big).
\end{align*}
By symmetry of~$K$, we have for any~$m$
\[
 \intl_{{\tilde F}^\omega(-\infty,a)}d\nuo(x)\intl_V d\nuo(y)
K^m(x,y) = \intl_V \Po^y[\xi_m\cdot\e < a]\, d\nuo(y),
\]
so Lemma~\ref{l_separated} now follows from~\eqref{bound_n_d>=4}.
\qed

Let us consider a sequence of i.i.d.\ random
variables~$Z_1,Z_2,Z_3,\ldots$ with uniform distribution on
$\{1,2,\ldots,N\}$ (where~$N$ is from Condition~T~(ii)), 
independent of
everything. Also, let us define ${\hat\xi}_n:=\xi_{Z_1+\cdots+Z_n}$.
Then, it is straightforward to obtain that, for
any~$x\in\partial\omega$ and~$B\subset \{y\in\partial\omega : -1\leq
(y-x)\cdot\e \leq 1\}$, we have
\begin{equation}
\label{oc_plotnost'}
 \Po^x[\hat\xi_1\in B] \geq N^{-1}r_1\nuo(B)
\end{equation}
for some $r_1>0$.
Let
\[
 {\hat K}(x,y) = \frac{1}{N}\sum_{j=1}^N K^j(x,y)
\]
be the transition density of the process $({\hat\xi}_n, n\geq 0)$.
Observe that this process is still reversible with the reversible
measure~$\nuo$, so that 
${\hat K}(x,y)={\hat K}(y,x)$ for all 
$x,y\in\partial\omega$.
Similarly to~\cite{CPSV2}, let us define
\begin{align}
 b(x) &= \Eo^x((\xi_1-x)\cdot\e)^2\nonumber\\
 &= \intl_{\partial\omega}((y-x)\cdot\e)^2
 K(x,y)\, \nuo(y),
\label{def_b}
\end{align}
and
\begin{align}
 {\hat b}(x) &= \Eo^x((\xi_{Z_1}-x)\cdot\e)^2\nonumber\\
 &= \intl_{\partial\omega}((y-x)\cdot\e)^2
 {\hat K}(x,y)\, \nuo(y).
\label{def_hat_b}
\end{align}
We suppose that $\tau(B)$ and~$\tau^+(B)$ are defined as 
in~\eqref{def_entrance}--\eqref{def_hitting} but with~$\xi$
instead of~$\bar\xi$, and let ${\hat\tau}(B)$ 
and~${\hat\tau}^+(B)$ be the corresponding quantities for the 
process~$\hat\xi$.
\begin{lmm}
\label{l_escape_upper}
 Suppose that $B,F\subset\partial\omega$ with
$\nuo(B)\in(0,\infty)$. 
Moreover, assume that $x\cdot\e\leq a$ for all $x\in B$, and
$y\cdot\e\geq a+u$ for all $y\in F$ (of course, the same result is
valid if we assume that $x\cdot\e\leq a$ for all $x\in F$, and
$y\cdot\e\geq a+u$ for all $y\in B$).
Then, there exist positive constants ${\tilde\gamma}_3$, 
${\tilde\gamma}_4$, such that
\begin{equation}
\label{eq_r_upper_xi}
 \nuo(B)\Po^B[\tau(F)<\tau^+(B)] \leq {\tilde\gamma}_3u^{-(d-1)}
   + \frac{1}{u^2}\intl_{{\tilde F}^\omega[a,a+u]}b(x)\, d\nuo(x),
\end{equation}
and
\begin{equation}
\label{eq_r_upper_hxi}
 \nuo(B)\Po^B[\hat\tau(F)<\hat\tau^+(B)] \leq
{\tilde\gamma}_4u^{-(d-1)}
   + \frac{1}{u^2}\intl_{{\tilde F}^\omega[a,a+u]}{\hat b}(x)\,
d\nuo(x).
\end{equation}
Moreover, \eqref{eq_r_upper_xi} and~\eqref{eq_r_upper_hxi} are valid
also in the finite tube~$\hD$ (in this case we assume that $a>0$
and $a+u<H$).
\end{lmm}

\noindent
\textit{Proof.}
We keep the notation~$\EEE(\cdot,\cdot)$ for the Dirichlet's form
with respect to~$K$, defined as in~\eqref{def_Dirichlet_form}. 
Suppose without restriction of generality that~$a=0$ and define the
function
\[
 h(x) = \begin{cases}
         0, &\text{ if }x\cdot\e \leq 0,\\
         1, &\text{ if }x\cdot\e \geq u,\\
         u^{-1}(x\cdot\e), &\text{ if }x\cdot\e \in (0,u).
        \end{cases}
\]
Using Proposition~\ref{p_Dirichlet} (observe that $h\in\mathcal{H}$)
and Lemma~\ref{l_separated}, we obtain
\begin{align*}
 \lefteqn{2\nuo(B) \Po^B[\tau(F)<\tau^+(B)]}~~~\\
 &\leq \EEE(h,h)\\
 &= \intl_{(\partial\omega)^2}d\nuo(x)\,d\nuo(y)K(x,y)
  (h(x)-h(y))^2\\
 &= 2\intl_{{\tilde F}^\omega(-\infty,0)}d\nuo(x)
  \intl_{{\tilde F}^\omega(u,\infty)}d\nuo(y) K(x,y)\\
 &\quad
   + u^{-2}\intl_{({\tilde F}^\omega[0,u])^2}
     d\nuo(x)\,d\nuo(y)K(x,y)((y-x)\cdot\e)^2\\
 &\quad
  + 2u^{-2}\intl_{{\tilde F}^\omega[0,u]} d\nuo(x)
 \intl_{{\tilde F}^\omega(-\infty,0)}d\nuo(y)
 K(x,y)(x\cdot\e)^2\\
 &\quad
 + 2\intl_{{\tilde F}^\omega[0,u]} d\nuo(x)
 \intl_{{\tilde F}^\omega(u,\infty)}d\nuo(y)
 K(x,y)(1-u^{-1}x\cdot\e)^2\\
 &\leq 2{\tilde\gamma}_3u^{-(d-1)}
  +2u^{-2}\intl_{{\tilde F}^\omega[0,u]}
     d\nuo(x)
     \intl_{\partial\omega}d\nuo(y)K(x,y)((y-x)\cdot\e)^2\\
&=2{\tilde\gamma}_3u^{-(d-1)}
   + \frac{1}{u^2}\intl_{{\tilde F}^\omega[a,a+u]}b(x)\, d\nuo(x),
\end{align*}
and this proves~\eqref{eq_r_upper_xi}. The proof
of~\eqref{eq_r_upper_hxi} is completely analogous.
\qed

\medskip

We now work in finite tube~$\hD$.
Let us use the abbreviations $U_n={\tilde F}^\omega[n-1,n)$,
and $V_n={\tilde F}^\omega[n,H)\cup\hDr$.
Observe that, by Condition~T~(i), we have that for some 
${\tilde\gamma}_5
\in (0,+\infty)$
\begin{equation}
\label{U_bounded}
\nuo(U_n) \geq {\tilde\gamma}_5
\end{equation}
for all~$n$ and for $\IP$-a.a.\ $\omega$.

To distinguish between the seconds moments of the projected jump
length in finite and infinite tubes, we modify our notations in the
following way.
For $x\in\partial\hD$, let $b_H(x)$ and ${\hat b}_H(x)$ be the
quantities defined as in~\eqref{def_b} and~\eqref{def_hat_b}, but in
the finite tube~$\hD$. Let us use the notations $b_\infty(x)$ and
${\hat b}_\infty(x)$ for the corresponding quantities in the
infinite tube.
Now, we need an estimate on the integrals
appearing in the right-hand sides of~\eqref{eq_r_upper_xi}
and~\eqref{eq_r_upper_hxi}, for the case of the finite tube:
\begin{lmm}
\label{l_int_b}
Suppose that $0<s_1<s_2<1$ and assume that $d\geq 3$ and Condition~T
holds. Then, we have
\begin{equation}
\label{eq_limsup}
 \limsup_{H\to\infty} \frac{1}{H}
      \intl_{{\tilde F}^\omega[s_1H,s_2H]}b_H(x)\, d\nuo(x) <\infty
 \qquad \IP\text{-a.s.},
\end{equation}
and the same is valid with~${\hat b}_H$ on the place of $b_H$.
\end{lmm}

\noindent
\textit{Proof.}
Let us recall some notations from~\cite{CPSV2}.
Define 
\[
 \gs = \{(\omega,u) : \omega\in\Omega, u\in\partial\omega_0\}.
\]
Define the probability measure~$\IQ$ on~$\gs$ by
\begin{equation}
\label{def_Q}
 d\IQ(\omega,u) = \frac{1}{\ZZ} \kappa_{0,u}^{-1} \, d\muo_{0}(u) \,
d\IP(\omega),
\end{equation}
where $\muo_0$ is the $(d-2)$-dimensional Hausdorff measure on the
boundary of~$\omega_0$,
$\kappa_{0,u}$ is the scalar product of the normal vectors pointing
inside the section and inside the tube (see Section~2
of~\cite{CPSV2} for details),
and~$\ZZ=\int_\Omega d\IP 
         \int_\Xi\kappa_{0,u}^{-1}d\muo_0(u)$ 
is the normalizing constant. In Lemma~3.1 of~\cite{CPSV2} it is
shown that~$\IQ$ is the invariant law
of the environment seen from the walker, that is
\begin{equation}
\label{IQ_is_invariant}
 \exq{\Eo[f(\theta_{\xi_n\cdot \e}\omega,\UU\xi_n)\mid
\xi_0=(0,u)]}=\exq{f}.
\end{equation}
Using also that
\[
 \Eo^x((\xi_n-x))^2 \leq n \sum_{k=1}^n \Eo^x((\xi_k-\xi_{k-1}))^2
\]
and~\eqref{IQ_is_invariant}, it is straightforward to obtain that
$\exq{b_\infty}<\infty$ implies $\exq{{\hat b}_\infty}<\infty$.
So, using the notations of~\cite{CPSV2}, by the ergodic theorem
we obtain
\begin{align}
 \frac{1}{H}
      \intl_{{\tilde F}^\omega[0,H]}b_\infty(x)\, d\nuo(x)
 &= \frac{1}{H} \intl_0^H d\alpha
  \intl_{\Xi}d\muo_\alpha(v)\kappa^{-1}_{\alpha,v}
 b_\infty(\theta_\alpha\omega,v)\nonumber\\
 &\to \exq{b_\infty} \qquad \text{as $H\to\infty$},
\label{conv_b}
\end{align}
a.s.\ and in~$L^1$,
and the same with ${\hat b}_\infty$ on the place of~$b_\infty$.
Then, \eqref{eq_limsup} follows from the fact that, for all~$H$,
$b_H(x)\leq b_\infty(x)$ for all $x\in \hF$. Now, 
with~${\hat b}_\infty$ instead
of~$b_\infty$, the previous inequality is not necessarily valid. So,
to prove~\eqref{eq_limsup} for~${\hat b}_H$ instead of~$b_H$,
consider~$x\in\partial\omega$ such that
$H^{-1}(x\cdot\e)\in[s_1,s_2]$, and write (note that for
all~$x\in\partial\hD$ we have ${\hat b}_H(x)\leq H^2$)
\begin{align*}
 \big|{\hat b}_H(x)-{\hat b}_\infty(x)\big| &\leq
 H^2 \Po^x\big[\max_{k\leq N}|(\xi_k-x)\cdot\e| \geq
(s_1\wedge(1-s_2))H^{-(d-1)}\big]\\
 &\leq C_1 H^{-(d-3)}
\end{align*}
(recall that $d\geq 3$), and then we obtain~\eqref{eq_limsup}
for~${\hat b}_H$ as well.
\qed

\medskip

Next, we obtain a lower bound for certain escape probabilities:
\begin{lmm}
\label{l_Dirichlet_lower}
 Suppose that $H/4\leq n\leq H-1$, and $m<n$. Also, assume that
$d\geq 3$ and Condition~T holds. Then, there
exist positive constants ${\tilde\gamma}_7$, ${\tilde\gamma}_8$, 
such that
\begin{equation}
\label{eq_low_U}
 \Po^{U_m}[\hat\tau(V_n)<{\hat\tau}^+(U_m)] \geq
\frac{{\tilde\gamma}_7}{n-m},
\end{equation}
and
\begin{equation}
\label{eq_low_D}
 \Po^{\hDl}[\hat\tau(V_n)<{\hat\tau}^+(\hDl)] \geq
\frac{{\tilde\gamma}_8}{H}.
\end{equation}
\end{lmm}

\noindent
\textit{Proof.}
Let~$\hat\EEE$ be the Dirichlet form corresponding to~$\hat K$
(cf.~\eqref{def_Dirichlet_form}).
First, let us prove~\eqref{eq_low_U}.
As in Proposition~\ref{p_Dirichlet}, we use the notation
${\hat h}(x)=\Po^x[\hat\tau(V_n)<{\hat\tau}(U_m)]$;
observe that ${\hat h}(x)=0$ for all $x\in U_m$
and ${\hat h}(y)=1$ for all $y\in V_n$ (and hence for all $y\in
U_{n+1}$). Using this fact together
with~\eqref{U_bounded} and Cauchy-Schwarz inequality, we write
(abbreviating $u:=n-m$)
\begin{align*}
 \lefteqn{2\nuo(U_m)\Po^{U_m}[\hat\tau(V_n)<{\hat\tau}^+(U_m)]}\\
 &= {\hat\EEE}({\hat h},{\hat h})\\
&\geq \sum_{j=0}^u \intl_{U_{m+j}}d\nuo(x_j)
\intl_{U_{m+j+1}}d\nuo(x_{j+1}) {\hat K}(x_j,x_{j+1})
  ({\hat h}(x_j)-{\hat h}(x_{j+1}))^2\\
&= \Big(\prod_{j=0}^{u+1}\nuo(U_{m+j})\Big)^{-1}
  \intl_{U_m}d\nuo(x_0)\ldots \intl_{U_{m+u+1}}d\nuo(x_{u+1})\\
&\qquad\qquad\qquad   \sum_{j=0}^{u}\nuo(U_{m+j})\nuo(U_{m+j+1})
    {\hat K}(x_j,x_{j+1})
  ({\hat h}(x_j)-{\hat h}(x_{j+1}))^2\\
&\geq 
 N^{-1}r_1{\tilde\gamma}_5^2
 \Big(\prod_{j=0}^{u+1}\nuo(U_{m+j})\Big)^{-1}
  \intl_{U_m}d\nuo(x_0)\ldots\\
&\qquad\qquad\qquad \ldots
\intl_{U_{m+u+1}}d\nuo(x_{u+1})\sum_{j=0}^{u}
  ({\hat h}(x_j)-{\hat h}(x_{j+1}))^2\\
&\geq \frac{N^{-1}r_1{\tilde\gamma}_5^2}{u+1}
  \Big(\prod_{j=0}^{u+1}\nuo(U_{m+j})\Big)^{-1}
\intl_{U_m}d\nuo(x_0)\ldots\intl_{U_{m+u+1}}d\nuo(x_{u+1})\\
&=  \frac{N^{-1}r_1{\tilde\gamma}_5^2}{n-m+1},
\end{align*}
and this proves~\eqref{eq_low_U}. 
By denoting ${\hat h}(x)=\Po^x[\hat\tau(V_n)<{\hat\tau}(\hDl)]$
and writing
\begin{align*}
 \lefteqn{2\nuo(\hDl)\Po^{\hDl}[\hat\tau(V_n)<{\hat\tau}^+(\hDl)]}\\
&\geq \sum_{j=1}^{n+1} \intl_{U_j}d\nuo(x_j)
\intl_{U_{j+1}}d\nuo(x_{j+1}) {\hat K}(x_j,x_{j+1})
  ({\hat h}(x_j)-{\hat h}(x_{j+1}))^2\\
 &\quad + \intl_{\hDl}d\nuo(x_0) \intl_{U_1}d\nuo(x_1) 
{\hat K}(x_0,x_1)
  ({\hat h}(x_0)-{\hat h}(x_1))^2
\end{align*}
in exactly the same way one can
show~\eqref{eq_low_D}. This concludes the proof of
Lemma~\ref{l_Dirichlet_lower}.
\qed

\medskip

Next, we need (pointwise) estimates on the probabilities of
exiting the tube at the left boundary:
\begin{lmm}
\label{l_gambler's_ruin}
Assume Condition~T and $d\geq 3$.
 Suppose also that $n\in (\frac{H}{4},\frac{3H}{4})$, and $m\in
(0,n]$. Then, there
exists ${\tilde\gamma}_9$ such that for all $x\in U_m$ we have
\begin{equation}
\label{gambler's_ruin}
 \Po^x[{\hat\tau}(\hDl)<{\hat\tau}(V_n)] 
      \leq \frac{{\tilde\gamma}_9(n-m+1)}{H}.
\end{equation}
\end{lmm}

\noindent
\textit{Proof.}
 From now on, we assume for technical reasons that $m>\frac{H}{8}$ 
(in any case, otherwise the upper bound~$1$ is good enough for us).
First, by Lemmas~\ref{l_escape_upper} and~\ref{l_int_b},  we obtain
that
\begin{equation}
\label{go_left}
 \Po^{U_m}[{\hat\tau}(\hDl)<{\hat\tau}^+(U_m)] \leq
\frac{C_1}{H}.
\end{equation}
Next, Lemma~\ref{l_Dirichlet_lower} implies that
\begin{equation}
\label{go_right}
 \Po^{U_m}[{\hat\tau}(V_n)<{\hat\tau}^+(U_m)] \geq
\frac{C_2}{n-m+1}.
\end{equation}
Also, from~\eqref{oc_plotnost'} it is clear that for any $x\in U_m$
we have
\begin{equation}
\label{podprygnuli}
 \Po^x[{\hat\tau}^+(U_m)<{\hat\tau}(\hDl\cup V_n)]
  \geq \Po^x[{\hat\xi}_1\in U_m] \geq C_3
\end{equation}
for some $C_3>0$. 

Now, denote $\sigma_0={\hat\tau}(U_m)$, 
$\sigma_{k+1}=\min\{j > \sigma_k : {\hat\xi}_j\in U_m\}$ to be
the successive times when the set~$U_m$ is visited. 
By Corollary~\ref{c_trans_dens}~(i) and~\eqref{podprygnuli},
we obtain that, conditional on not hitting $\hDl\cup V_n$, the
process of successive returns to~$U_m$ is reversible with the
reversible density $\pi_m(x)$, such that for all $x\in U_m$
\[
 C_4 \leq \pi_m(x) \leq C_5
\]
for some positive constants $C_4,C_5$.
Using also~\eqref{go_left} and~\eqref{go_right}, we obtain that
there are constants $C_6,C_7>0$ such that for any~$k$
\begin{align*}
 \Po^{U_m}[{\hat\tau}(\hDl)<{\hat\tau}^+(U_m)
         \mid {\hat\tau}(\hDl\cup V_n)>\sigma_k] &\leq
\frac{C_6}{H},\\
 \Po^{U_m}[{\hat\tau}(V_n)<{\hat\tau}^+(U_m)
         \mid {\hat\tau}(\hDl\cup V_n)>\sigma_k] &\geq
\frac{C_7}{n-m+1}.
\end{align*}
So, we can write
\begin{align}
 \Po^{U_m}[{\hat\tau}(\hDl)<{\hat\tau}(V_n)] 
  &= \sum_{k=1}^\infty 
   \Po^{U_m}\big[{\hat\tau}(\hDl)<{\hat\tau}(V_n) 
      \mid {\hat\tau}(\hDl\cup V_n)\in(\sigma_{k-1},\sigma_k]\big]
\nonumber\\
  &\qquad\qquad \times
\Po^{U_m}\big[{\hat\tau}(\hDl\cup V_n)
           \in(\sigma_{k-1},\sigma_k]\big]
\nonumber\\
 &\leq \sum_{k=1}^\infty \frac{C_6/H}{C_7/(n-m+1)}
\Po^{U_m}\big[{\hat\tau}(\hDl\cup V_n)
           \in(\sigma_{k-1},\sigma_k]\big]
\nonumber\\
&= \frac{C_6C_7^{-1}(n-m+1)}{H}.
\label{bound_set}
\end{align}

Now, the ``pointwise'' version of~\eqref{bound_set} is substantially
more difficult to prove.

Consider a sequence of i.i.d.\ random variables $\zeta_n\in\{0,1\}$
with 
\[
 P[\zeta_n=1] = N^{-1}r_1{\tilde\gamma}_5
\]
(recall~\eqref{oc_plotnost'} and~\eqref{U_bounded}). Then,
one can couple the random sequences
$({\hat\xi}_n, n\geq 1)$ with $\zeta=(\zeta_n, n\geq 1)$ in
such a
way that when the event $\{\zeta_n=1\}$ occurs, ${\hat\xi}_n$
has the stationary distribution on $U_{[{\hat\xi}_{n-1}\cdot\e]}$.
We denote by~$\Pozz$ and~$\Eozz$ the probability and expectation
with fixed~$\omega$ and~$\zeta$, and let~$E^{\zeta}$ be the
expectation with respect to~$\zeta$. 
One can
formally define~$\Pozz$ in the following way. 
For any $x\in U_i$, define the transition density~$R_x$ by
\[
 (1-N^{-1}r_1{\tilde\gamma}_5) R_x(y) = 
  \begin{cases}
   K(x,y),& \text{if }   y\notin U_i,\\
   K(x,y)-\frac{N^{-1}r_1{\tilde\gamma}_5}{\nuo(U_i)},
    & \text{if }   y\in U_i.
  \end{cases}
\]
Let $\mathcal{R}_x$ be the distribution on~$\partial\omega$ with the
density~$R_x$, and let $\mathcal{U}_i$ be the uniform distribution
on~$U_i$. Then, given ${\hat\xi}_{n-1}=x\in U_i$, the law of
${\hat\xi}_n$ under~$\Pozz$ is given by
\[
 \1{\zeta_n=1}\mathcal{U}_i + \1{\zeta_n=0}\mathcal{R}_x.
\]
Also, let us define ${\hat\kappa}=\min\{n\geq 1: \zeta_n=1\}$.

Now, observe that
\begin{equation}
\label{=j}
 [{\hat\xi}_j\cdot\e] = [{\hat\xi}_{j-1}\cdot\e]
\text{ on } \{j={\hat\kappa}\}
\end{equation}
and, for~$i$ such that $i<j$,
\begin{align}
 \lefteqn{E^{\zeta} \big(\Pozz^x\big[|({\hat\xi}_i-
  {\hat\xi}_{i-1})\cdot\e| \geq u \big]  \mid {\hat\kappa}=j
\big)}
\nonumber\\
 &= E^{\zeta} \big(\Pozz^x\big[|({\hat\xi}_i-
  {\hat\xi}_{i-1})\cdot\e| \geq u
\big]  \mid \zeta_i=0 \big)
\nonumber\\
&\leq \frac{1}{P^{\zeta}[\zeta_i=0]}
  \Po^x\big[|({\hat\xi}_i-
  {\hat\xi}_{i-1})\cdot\e| \geq u \big]
\nonumber\\
&\leq C_8h^{-(d-1)},
\label{<j}
\end{align}
recall~\eqref{bound_n_d>=4}. 
Then, write using~\eqref{=j} and~\eqref{<j}
\begin{align}
\lefteqn{E^{\zeta} \Pozz^x\big[ \max_{\ell\leq{\hat\kappa}}
|({\hat\xi}_\ell-{\hat\xi}_0)\cdot\e|\geq s\big]}
\nonumber\\
 &= \sum_{j=1}^\infty P^{\zeta}[{\hat\kappa}=j]
  E^{\zeta}\big(\Pozz^{U_0}\big[ \max_{\ell\leq{\hat\kappa}}
 |({\hat\xi}_\ell-{\hat\xi}_0)\cdot\e|
\geq s\big]\mid {\hat\kappa}=j\big)
\nonumber\\
&\leq \sum_{j=1}^\infty P^{\zeta}[{\hat\kappa}=j]
  E^{\zeta}\big(\Pozz^x\big[\text{there exists }i\leq j
\text{ such that }
\nonumber\\
& \qquad\qquad\qquad\qquad\qquad\qquad\qquad
|({\hat\xi}_i-{\hat\xi}_{i-1})\cdot\e|
\geq s/j\big]\mid {\hat\kappa}=j\big)
\nonumber\\
&\leq \sum_{j=1}^\infty P^{\zeta}[{\hat\kappa}=j]
jC_9\Big(\frac{s}{j}\Big)^{-(d-1)}
\nonumber\\
&= C_9 s^{-(d-1)} \sum_{j=1}^\infty j^d P^{\zeta}[{\hat\kappa}=j]
\nonumber\\
&= C_{10} s^{-(d-1)}.
\label{oc_kappa_1}
\end{align}
Now, using~\eqref{oc_kappa_1}, we have for an arbitrary~$x\in U_m$
\begin{align}
 \lefteqn{\Po^x[{\hat\tau}(\hDl)<{\hat\tau}(V_n)]}
 \nonumber\\
  &= E^{\zeta}\Pozz^x[{\hat\tau}(\hDl)<{\hat\tau}(V_n)]
 \nonumber\\
&\leq E^{\zeta} \Pozz^x\big[\max_{j\leq {\hat\kappa}} 
|(x-{\hat\xi}_j)\cdot\e|<H/16,
{\hat\tau}(\hDl)<{\hat\tau}(V_n)\big]
 \nonumber\\
&\qquad + E^{\zeta} \Pozz^x\big[\max_{j\leq {\hat\kappa}} 
|(x-{\hat\xi}_j)\cdot\e|\geq H/16\big]
\nonumber\\
&\leq E^{\zeta} \Pozz^x\big[\max_{j\leq {\hat\kappa}} 
|(x-{\hat\xi}_j)\cdot\e|<H/16,
{\hat\tau}(\hDl)<{\hat\tau}(V_n)\big] + C_{11}H^{-(d-1)}.
\label{I+II}
\end{align}
Let us deal with the first term in~\eqref{I+II}.
We have, taking advantage of~\eqref{bound_set}
and~\eqref{oc_kappa_1} (recall that $d\geq 3$)
\begin{align*}
 \lefteqn{E^{\zeta} \Pozz^x\big[\max_{j\leq {\hat\kappa}} 
|(x-{\hat\xi}_j)\cdot\e|<H/16,
{\hat\tau}(\hDl)<{\hat\tau}(V_n)\big]}\\
&\leq
\sum_{\ell\geq H/16}
E^{\zeta}\Pozz^x\big[[{\hat\xi}_{{\hat\kappa}}]=\ell\big]
\Po^{U_\ell}[{\hat\tau}(\hDl)<{\hat\tau}(V_n)] \\
&\leq \frac{C_{12}(n-m+1)}{H}\sum_{\ell\geq m}
E^{\zeta}\Pozz^x\big[[{\hat\xi}_{{\hat\kappa}}]=\ell\big]\\
&\quad +
\sum_{\frac{H}{16}\leq\ell<m}
\frac{C_{12}\big((n-m+1)+(m-\ell)\big)}{H}
E^{\zeta}\Pozz^x\big[[{\hat\xi}_{{\hat\kappa}}]=\ell\big]\\
&\leq \frac{C_{13}(n-m+1)}{H},
\end{align*}
and this concludes the proof of Lemma~\ref{l_gambler's_ruin}.
\qed

\medskip

Next, we prove a result which shows that it is unlikely that a
particle crosses the tube~$\hD$ ``too quickly''. 
Suppose that one particle is injected (uniformly) at random
at~$\hDl$ into the tube~$\hD$, and we still denote by~$\CC_H$ 
the event that
it crosses the tube without going back to~$\hDl$, i.e.,
$\CC_H=\{\tau(\hDr)<\tau^+(\hDl)\}$ (one can see that there is no
conflict with the notation of Section~\ref{s_perm}). Also,
recall that~$\TT_H$ stands for the total
lifetime of the particle as defined in Section~\ref{s_perm}, i.e.,
if~$X_t$ is the location of the
particle at time~$t$, then $\TT_H=\min\{t>0:X_t\in \hDl\cup\hDr\}$.
\begin{lmm}
\label{l_cross}
 For any~$\eps>0$ there exists (large enough)~$m$ with the following
property: there exists large enough~$H_0=H_0(\omega)$ such that for
all $H\geq H_0$
\begin{equation}
\label{eq_cross}
 \Po^{\hDl}[\CC_H, \TT_H\leq m^{-1}H^2] \leq \frac{\eps}{H}.
\end{equation}
\end{lmm}

\noindent
\textit{Proof.}
For $H,m,\eps_1>0$, we say that $x\in\partial\omega$ is
$(H,m,\eps_1)$-good if
\begin{equation}
\label{def_Hme_good} 
 \Po^x\big[\sup_{t\leq m^{-1}H^2}|(X_t-x)\cdot\e|<H/4\big]
           \geq 1-\eps_1.
\end{equation}
 Let~$L\in\Z$ be a large positive parameter to be
specified later; for $n\in\Z$ denote 
$I_n={\tilde F}^\omega[nL,(n+1)L)$; denote also
\[
{\tilde I}_n^{\eps_1} = \{x\in I_n: x\text{ is not
$(H,m,\eps_1)$-good}\}.
\]
Now, consider first the case $d\geq 3$. From now on we suppose
that~$m$
is sufficiently large to assure the following: 
\[
 P\big[\sup_{t\leq m^{-1}}|B_t|<1/4\big] \geq 1-\frac{\eps_1}{2},
\]
where $B_t$ is the standard Brownian motion and~$P$
is the corresponding probability measure. 
In this case, if the invariance
principle holds, then for any fixed~$\eps_1>0$ every~$x$ is
$(H,m,\eps_1)$-good for all large enough~$H$. Using the monotone
convergence theorem, it is straightforward to obtain that
for fixed $L,m,\eps_1,\eps_2$ there exists large enough~$H_0$ such
that for all $H\geq H_0$
\begin{equation}
\label{I_good}
 \IP[\nuo({\tilde I}_0^{\eps_1})<\eps_2] > \frac{3}{4}.
\end{equation}
Then, by the ergodic theorem, there exists large enough~$H_0$ such
that for all $H\geq H_0$ there exists $n_0=n_0(H)$ such that
$I_{n_0}\subset {\tilde F}^\omega(H/4,3H/4)$, and 
$\nuo({\tilde I}_{n_0}^{\eps_1})<\eps_2$.

Now, let us consider also the
event~${\hat\CC}_H=\{{\hat\tau}(V_{Ln_0})<{\hat\tau}^+(\hDl)\}$ (that
is, with respect to the process~${\hat\xi}$, the particle
enters~$V_{Ln_0}$ before coming back to~$\hDl$).
Then, write
\begin{align}
 \Po^{\hDl}[\CC_H, \TT_H\leq m^{-1}H^2] &\leq
 \Po^{\hDl}[{\hat\CC}_H, \TT_H\leq m^{-1}H^2] 
+ \Po^{\hDl}[{\hat\CC}_H^c,\CC_H]
\nonumber\\
 &= \Po^{\hDl}[{\hat\CC}_H]\Po^{\hDl}[\TT_H\leq m^{-1}H^2\mid
{\hat\CC}_H] \nonumber\\
&\quad + \Po^{\hDl}[\CC_H]\Po^{\hDl}[{\hat\CC}_H^c\mid\CC_H].
\label{cross_decompose}
\end{align}
Now, by Lemmas~\ref{l_escape_upper} and~\ref{l_int_b}, we can write
for some~$C_1>0$
\begin{equation}
\label{oc_peresech'}
 \max\{\Po^{\hDl}[\CC_H],\Po^{\hDl}[{\hat\CC}_H]\} \leq
\frac{C_1}{H}.
\end{equation}
Then, from~\eqref{bound_n_d>=4} we obtain that
\begin{equation}
\label{jump_back}
\Po^{\hDl}[{\hat\CC}_H^c\mid\CC_H] \leq 
\sup_{x\in\hDr}
\Po^x\big[\max_{j\leq N}|\xi_j\cdot\e-H|\geq H/4\big]
 \leq C_2 H^{-(d-1)}
\end{equation}
for some $C_2>0$. So, to complete the proof of~\eqref{eq_cross}, it
remains to prove that the term
$\Po^{\hDl}[\TT_H\leq m^{-1}H^2\mid {\hat\CC}_H]$
in~\eqref{cross_decompose} is small.

To do this, let us recall that, by Lemma~\ref{l_hitting_revers}~(i),
for any~$F'\subset V_{Ln_0}$, we have
\begin{equation}
\label{come_to_F'}
  \Po^{\hDl}[\xi_{{\hat\tau}(V_{Ln_0})}\in F' \mid {\hat\CC}_H]
 = \big(\nuo(\hDl)\Po^{\hDl}[\hat\CC_H]\big)^{-1}
 \intl_{F'}\Po^y[{\hat\tau}(\hDl)<{\hat\tau}^+(V_{Ln_0})]\,d\nuo(y).
\end{equation}
By Lemma~\ref{l_Dirichlet_lower}, we have that for some $C_3>0$
\begin{equation}
\label{oc_nu_P}
 \big(\nuo(\hDl)\Po^{\hDl}[\hat\CC_H]\big)^{-1} \leq C_3 H.
\end{equation}
For $j\geq 1$ denote $S_j=\hDl\cup U_1\cup\ldots\cup U_j$. 
Using Lemma~\ref{l_gambler's_ruin}, we can write 
for any~$y\in V_{Ln_0}$
\begin{align}
 \Po^y[{\hat\tau}(\hDl)<{\hat\tau}^+(V_{Ln_0})] &=
  \intl_{\partial\omega}{\hat K}(y,z)
  \Po^z[{\hat\tau}(\hDl)<{\hat\tau}(V_{Ln_0})]\,d\nuo(z)
\nonumber\\
 &\leq \intl_{\hDl}{\hat K}(y,z)\,d\nuo(z) \nonumber\\
 &\qquad  + \sum_{j=1}^{Ln_0} \frac{{\tilde\gamma}_9(Ln_0-j+1)}{H}
  \intl_{U_j}{\hat K}(y,z)\,d\nuo(z)\nonumber\\
 & \leq \frac{{\tilde\gamma}_9}{H}\sum_{j=1}^{Ln_0}
    \intl_{S_j}{\hat K}(y,z)\,d\nuo(z).
\label{upper_Py}
\end{align}
So, by~\eqref{bound_n_d>=4}, in the case $d\geq 3$, we obtain
 from~\eqref{upper_Py} that for some positive constant~$C_4$
\[
 \Po^y[{\hat\tau}(\hDl)<{\hat\tau}^+(V_{Ln_0})] \leq \frac{C_4}{H}
\]
and, by~\eqref{come_to_F'}, \eqref{oc_nu_P}, and the construction
of~$n_0$ we obtain that
\begin{equation}
\label{oc_tilde_I}
\Po^{\hDl}[\xi_{{\hat\tau}(V_{Ln_0})}\in {\tilde I}_{n_0}^{\eps_1}
   \mid {\hat\CC}_H] \leq C_3C_4\eps_2.
\end{equation}
Next, integrating~\eqref{upper_Py} over $V_{n_0}\setminus I_{n_0}$,
we
obtain from Lemma~\ref{l_separated} that
\begin{align*}
 \intl_{V_{n_0}\setminus I_{n_0}}
\Po^y[{\hat\tau}(\hDl)<{\hat\tau}^+(V_{Ln_0})]\,d\nuo(y)
 &\leq \frac{{\tilde\gamma}_9}{H}\sum_{j=1}^{Ln_0}
    \intl_{V_{n_0}\setminus I_{n_0}}d\nuo(y)
 \intl_{S_j}d\nuo(z){\hat K}(y,z)\\
&\leq \frac{C_5{\tilde\gamma}_9}{H}
  \sum_{j=1}^{Ln_0} (Ln_0+L-j)^{-(d-1)}\\
&\leq \frac{C_6}{H} L^{-(d-2)}.
\end{align*}
Again using~\eqref{come_to_F'}, \eqref{oc_nu_P}, we obtain that
\begin{equation}
\label{oc_V-I}
\Po^{\hDl}[\xi_{{\hat\tau}(V_{Ln_0})}\in V_{n_0}\setminus I_{n_0}
\mid {\hat\CC}_H] \leq C_3C_6 L^{-(d-2)}.
\end{equation}
So, \eqref{oc_tilde_I} and~\eqref{oc_V-I} imply that for
any~$\eps_3>0$ there exists large enough~$L$
such that for all large enough~$H$ we have
\[
 \Po^{\hDl}[\xi_{{\hat\tau}(V_{Ln_0})}\in I_{n_0}\setminus 
{\tilde I}_{n_0}^{\eps_1} \mid {\hat\CC}_H] 
\geq 1-\eps_3.
\]
But then, since all 
$x\in I_{n_0}\setminus {\tilde I}_{n_0}^{\eps_1}$ are
$(H,m,\eps_1)$-good, from~\eqref{def_Hme_good} we obtain that
\begin{equation}
\label{oc_4th_term}
 \Po^{\hDl}[\TT_H\leq m^{-1}H^2\mid {\hat\CC}_H] \leq 
   1-(1-\eps_1)(1-\eps_3).
\end{equation}
Using~\eqref{oc_peresech'}, \eqref{jump_back},
and~\eqref{oc_4th_term} in~\eqref{cross_decompose}, we conclude the
proof of~\eqref{eq_cross} in the case $d\geq 3$.

Let us prove the lemma in the case $d=2$. Take
\[
 L=\sup\{|(x-y)\cdot\e| : x,y\in\RR, x\near y\}.
\]
Note that $b_H(x)\leq L^2$, so Lemma~\ref{l_escape_upper} implies
that $\Po^{\hDl}[\CC_H]\leq C_7H^{-1}$ for some $C_7>0$.
By Condition~T~(iii) we obtain that
\[
\Po^{\hDl}[\xi_{{\hat\tau}(V_{Ln_0})}\in V_{n_0}\setminus I_{n_0}
\mid \CC_H] = 0,
\]
and, since for any $x\in I_{Ln_0}$, 
$y\in{\tilde F}^\omega[0,L(n_0-1))$ we have $K(x,y)=0$, 
we then obtain
\[
 \Po^{\hDl}[\xi_{{\hat\tau}(V_{Ln_0})}\in I_{n_0}\setminus 
{\tilde I}_{n_0}^{\eps_1} \mid \CC_H] 
\geq 1-\eps_4
\]
for a small~$\eps_4>0$. The
proof of~\eqref{eq_cross} in the case $d=2$ then follows in the same
way.
\qed

\section{On the steady state of the Knudsen gas}
\label{s_pr_steady}
In this section we prove the theorem that characterizes the
stationary regime for the Knudsen gas in a finite tube.

\medskip
\noindent
\textit{Proof of Theorem~\ref{t_stat_measure}.}
In order to prove item~(i), we consider the process with
absorbing/injection boundaries in \emph{both} $\hDl$ and~$\hDr$
(that is, the injection is given by two independent
Poisson processes in $\hDl \times \Sph_\e$ and $\hDr \times
\Sph_{(-\e)}$ with intensities $|\Sph^{d-1}|^{-1} \lambda |\e \cdot
u|\, dx\, du$ in both cases).

Fix a sequence of positive numbers $u_k\nearrow \infty$ such that
$\lambda u_k\in\Z$ for all~$k$.
For each~$k$, consider a domain~$\Phi_k$ with the following
properties
\begin{itemize}
 \item $\hD\subset \Phi_k$, $\hF\subset\partial\Phi_k$,
$(\hDl\cup\hDr)\subset\Phi_k$;
 \item $|\Phi_k|=u_k$;
 \item any segment $ab$ with $a\in \hDl\cup\hDr$, $b\in
\partial\Phi_k\setminus\hF$
 has length at least $u_k^{1/(2d)}$
\end{itemize}
\begin{figure}
 \centering
\includegraphics{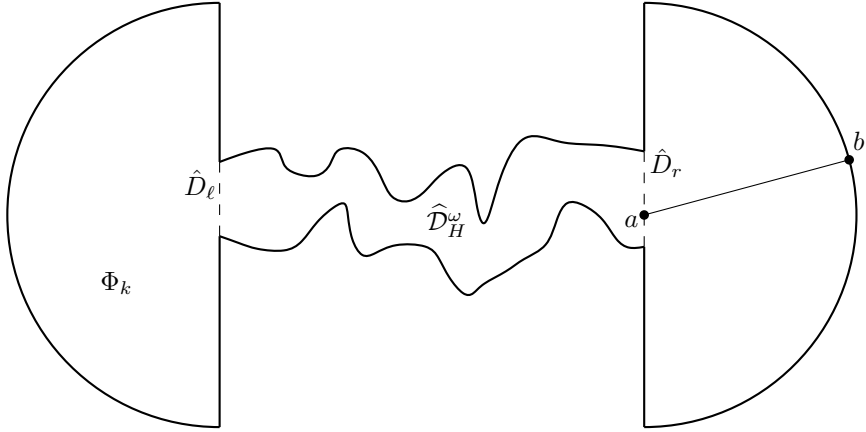}
\caption{On the construction of domain $\Phi_k$}
\label{f_reservoir}
\end{figure}
(one may construct such a domain e.g.\ as shown on
Figure~\ref{f_reservoir}).
Now, let us consider $\lambda u_k$ independent particles in~$\Phi_k$.
By Theorem~2.4 of~\cite{CPSV1}, the unique invariant measure of this
system
is product of uniform measures in location and direction.
We are going to compare this process (observed only
on~$\hD$) with the process with absorbing/injection boundaries in 
both $\hDl$ and~$\hDr$ (naturally, we assume that the injection
is with the cosine law and with the same intensity mentioned in
Theorem~\ref{t_stat_measure}.
Let~$E^{(k)}$ be the expectation for the above process in~$\Phi_k$
with~$\lambda u_k$ particles, with respect to the invariant measure. 
Also, we denote by~$E_t$ the expectation with respect to the 
process with absorbing/injection boundaries in $\hDl\cup\hDr$ at
time~$t$, with the initial configuration chosen from the Poisson
point process in $\hD\times \Sph^{d-1}$
with intensity~$\lambda|\Sph^{d-1}|^{-1}$.

Let~$\psi$ be a function on~$\hD\times\Sph^{d-1}$, taking values on
the interval~$[0,1]$. For a configuration
$\eta=(x_1,v_1,\ldots,x_r,v_r)$ in $\hD\times\Sph^{d-1}$
 (which means that we have~$r$ particles with positions
$x_1,\ldots,x_r\in\hD$ and vector speeds
$v_1,\ldots,v_r\in\Sph^{d-1}$), write
\[
 \psi(\eta) = \prod_{j=1}^r \psi(x_j,v_j).
\]
Denote also by
\[
 \bar{\psi} = \frac{1}{|\hD||\Sph^{d-1}|}
\intl_{\hD\times\Sph^{d-1}}\psi(x,v)\, dx\,dv
\]
the mean value of~$\psi$ on $\hD\times\Sph^{d-1}$.

Clearly, we have
\begin{align}
 E_0 \psi(\eta) &= e^{-\lambda|\hD|}\sum_{j=0}^\infty 
          \frac{(\lambda|\hD|)^j}{j!}\bar{\psi}^j\nonumber\\
  &= \exp\big(\lambda|\hD|(\bar{\psi}-1)\big) .
\label{gen_functional}
\end{align}
Also, it is straightforward to obtain that
\begin{equation}
\label{h_bino}
 E^{(k)}\psi(\eta) = \sum_{j=0}^{\lambda u_k} \binom{\lambda
u_k}{j}\Big(\frac{|\hD|}{u_k}\Big)^j
        \Big(1-\frac{|\hD|}{u_k}\Big)^{\lambda u_k-j} \bar{\psi}^j.
\end{equation}
Since, as $k\to\infty$, the binomial distribution with parameters
$\lambda u_k$ and $|\hD|/u_k$
converges to the Poisson distribution with parameter $\lambda|\hD|$,
for any~$\psi$ we have
\begin{equation}
\label{h_Poisson}
 \lim_{k\to\infty}E^{(k)}\psi(\eta) = E_0 \psi(\eta).
\end{equation}

Now, let us fix~$t_0$ and prove that for any~$\eps>0$
\begin{equation}
\label{h_time_t}
 \big|E^{(k)}\psi(\eta) - E_{t_0}\psi(\eta)\big| < \eps
\end{equation}
for all large enough~$k$.
For this, denote by $N^{(r)}(t_0)$ the total number of particles
which entered~$\hD$ through the right boundary~$\hDr$ up to
time~$t_0$.
For the process with absorption/injection, an elementary calculation
shows that $N^{(r)}(t_0)$ has Poisson
distribution with parameter $(\gamma_d|\Sph^{d-1}|)^{-1}\lambda
t_0|\hDr|$.
Let us suppose without restriction of generality that
$t_0<u_k^{1/(2d)}$ and denote
\[
 \Theta(t_0) = \{(x,v)\in\R^d\times\Sph_{(-\e)}: \text{ there exists
}t\in[0,t_0]
  \text{ such that }x+vt\in\hDr\};
\]
observe that $\Theta(t_0)\subset\Phi_k\times\Sph_{(-\e)}$. 

Now, a particle starting in~$x\in\Phi_k\setminus\hD$ with the
direction~$v$ will
cross~$\hDr$ by time~$t_0$ iff $(x,v)\in\Theta(t_0)$. So, 
it is straightforward to obtain that, for the process
in~$\Phi_k$, the random variable $N^{(r)}(t_0)$ has the binomial
distribution with parameters~$\lambda u_k$ and $\frac{t_0
|\hDr|}{\gamma_d |\Sph^{d-1}|u_k}$,
which converges to the Poisson distribution with parameter
$(\gamma_d|\Sph^{d-1}|)^{-1}\lambda t_0|\hDr|$
as $k\to\infty$. Then, conditioned on $\{N^{(r)}(t_0) = n\}$, for
both processes the~$n$ entering particles to~$\hDr$ (seen as a point
process on $\hDr\times\Sph_{(-\e)}\times[0,t_0]$)
are independent, each having
density $f(x,v,t)=t_0^{-1}|\hDr|^{-1}\gamma_d |v\cdot\e|$. 
Observe that the same considerations apply also to the particles
which enter through~$\hDl$. To obtain~\eqref{h_time_t}, we use now 
the following coupling argument. First of all, as we already know,
the initial configurations restricted to~$\hD$ 
for both processes can
be successfully coupled with probability
that converges to~$1$ as $k\to\infty$. Then, by the argument we just
presented, the same applies for the process of particles entering
through $\hDl\cup\hDr$. This shows that, with large probability,
both processes can be successfully coupled.

Now, combining~\eqref{h_Poisson} with~\eqref{h_time_t} and using the
fact that a point process is uniquely determined by its
characteristic functional (cf.\ e.g.\ Section~5.5 of~\cite{DV}), we
obtain that the Poisson point process in $\hD\times \Sph^{d-1}$
with intensity~$\lambda|\Sph^{d-1}|^{-1}$ is invariant for the
Knudsen gas with absorption/injection in $\hDr\cup\hDl$.

As for the convergence to the stationary state and the uniqueness,
this follows from an easy coupling argument.
Indeed, consider one process starting from the invariant measure
defined above, and another process starting from 
an arbitrary (fixed)
configuration. The initial particles are independent, but the newly
injected particles are the same for both processes. Then, since any
fixed particle will eventually disappear, the coupling time is a.s.\
finite, and so the system converges to the unique stationary state.
(Using Theorem~2.1 of~\cite{CPSV1}, with some more work one can show
that, for \emph{fixed} tube, this convergence is exponentially fast;
however, we do not need this kind of result in the present paper.)
This concludes the proof of the part~(i).

Let us prove the part~(ii).
Still considering the process with absorption and injection in
$\hDr\cup\hDl$, suppose that the particles entering through $\hDr$
are coloured red, and the particles
entering through~$\hDl$ 
are coloured green. So, we need to compute the
stationary measure for green particles. Using the (quasi)
reversibility of Knudsen stochastic billiard
(see Theorem~2.5 of~\cite{CPSV1}), we obtain that, given that there
is a particle in~$x=(\alpha,u)$ with the vector speed~$h$, the
probability that it is green equals
\[
 \Po^{(\alpha,u),-h}[\wp_{-\alpha}(X\cdot\e)<\wp_{H-\alpha}
            (X\cdot\e)].
\]
Using also the part~(i), we obtain that,
for the gas with injection only in~$\hDl$, the stationary measure is
that of Poisson point process with intensity
\[
\lambda |\Sph^{d-1}|^{-1} \Po^{(\alpha,u),-h}[\wp_{-\alpha}(X\cdot
\e)<\wp_{H-\alpha}(X\cdot\e)] \,d\alpha\, du\, dh.
\]
Note also that convergence and uniqueness follow from the same
coupling argument as in part~(i).
This concludes the proof of Theorem~\ref{t_stat_measure}.
\qed

Let us observe also that Theorem~\ref{t_stat_measure} 
allows us to characterize the stationary measure for Knudsen gas
where the injection takes place
from both sides, but with different intensities (which are constant
on~$\hDl$ and~$\hDr$).
We have
\begin{cor}
\label{cor_different}
 Consider now Knudsen gas with injection from both sides, with
respective intensities
$(\gamma_d|\Sph^{d-1}|)^{-1}\lambda$ and
$(\gamma_d|\Sph^{d-1}|)^{-1}\mu$ on~$\hDl$ and~$\hDr$
(without restriction of generality, let us suppose that 
$\lambda\geq \mu$). Then,
a Poisson point process with intensity measure
\[
 |\Sph^{d-1}|^{-1}\big(\mu+ 
(\lambda-\mu)  \Po^{(\alpha,u),-h}[\wp_{-\alpha}(X\cdot
\e)<\wp_{H-\alpha}(X\cdot\e)]\big)
  \,d\alpha\, du\, dh
\]
is the steady state of the Knudsen gas.
\end{cor}

\medskip
\noindent
\textit{Proof.} Indeed, one may imagine that particles of type~$1$
are injected from both sides
with intensity~$(\gamma_d|\Sph^{d-1}|)^{-1}\mu$ and particles of
type~$2$ are injected only from the left with
intensity~$(\gamma_d|\Sph^{d-1}|)^{-1}(\lambda-\mu)$,
and use Theorem~\ref{t_stat_measure}.
\qed

\section{Proofs of the results on transport diffusion and crossing
time}
\label{s_proofs_final}

\subsection{Proof of Theorems~\ref{t_dens_grad} and~\ref{t_current}}
\label{s_pr_trans}
For integers $i,j,\ell \geq 0$ define 
\begin{align*}
 R_{i,j}(g) &= \1{\wp_{-i}(g)>\wp_j(g)},\\
 G_{i,j,\ell}(g) &= \1{\wp_i(g)\leq \ell, \wp_i(g)<\wp_{-j}(g)}.
\end{align*}
Let~$B^{({\hat\sigma})}$ be a Brownian motion with diffusion
constant~${\hat\sigma}$, starting from the origin; we define
(being~$E$ the expectation with respect to the probability measure
on the space where the Brownian motion is defined)
\begin{align*}
 {\tilde R}_{i,j} &= E R_{i,j}(B^{({\hat\sigma})}) = \frac{i}{i+j},\\
 {\tilde G}_{i,j,\ell} &= E G_{i,j,\ell}(B^{({\hat\sigma})})
\end{align*}
to be the probabilities of the corresponding events for this Brownian
motion.

Fix an integer~$m$. For $(z,h)\in\omega\times\Sph^{d-1}$ and
$\eps_1>0$
define
\begin{align}
 T_\omega^{\eps_1}(z,h) =&\inf\big\{ s_0\geq 0 :
|\E^{z,h}R_{i,j}({\hat Z}^{(s)})
 -{\tilde R}_{i,j}|<\eps_1, \nonumber\\
 &\qquad\quad |\E^{z,h}G_{i,j,m}({\hat Z}^{(s)})
 -{\tilde G}_{i,j,m}|<\eps_1, \nonumber\\
 &\qquad\quad \text{ for all } i,j>0 \text{ such that } i+j=m, 
  \text{ and all }s\geq s_0\big\}. \label{def_Teps}
\end{align}
Intuitively, $T_\omega^{\eps_1}(z,h)$ is the scaling factor
one needs to use in order to assure that the rescaled (and projected
on~$e$)
trajectory of the Knudsen stochastic billiard stays 
sufficiently close to the Brownian motion.

By the portmanteau theorem, observe that, 
if the Knudsen stochastic billiard starting from $(z,h)$ satisfies
the quenched
invariance principle, this means that for any $\eps_1>0$ it holds
that $T_\omega^{\eps_1}(z,h)<\infty$.
Since, for $\IP$-almost every~$\omega$, the invariance principle
holds for a.a.\ starting
points $(z,h)$, we have
\[
 \intl_{\Omega} d\IP \;
|\Sph^{d-1}|^{-1}|\omega_0|^{-1}\intl_{\omega_0}du\intl_{\Sph^{d-1}}
dh\, 
 \1{T_\omega^{\eps_1}\big((0,u),h\big)<\infty} = 1.
\]
By the monotone convergence theorem, we obtain that
 for all $\eps_1,\eps_2>0$ there exists $t_{\eps_1,\eps_2}$ such
that 
\begin{equation}
\label{prob_good_point}
 \intl_{\Omega} d\IP \;
|\Sph^{d-1}|^{-1}|\omega_0|^{-1}\intl_{\omega_0}du\intl_{\Sph^{d-1}}
dh\, 
 \1{T_\omega^{\eps_1}\big((0,u),h\big)\leq t_{\eps_1,\eps_2}} \geq
1-\eps_2.
\end{equation}
So, using the Ergodic Theorem, we obtain for almost all~$\omega$ and
all~$H$ large enough 
\begin{equation}
\label{good_points}
 |\Sph^{d-1}|^{-1} \big|\{(z,h)\in \hD\times \Sph^{d-1} :
   T_\omega^{\eps_1}(z,h)>t_{\eps_1,\eps_2}\}\big| \leq 2\eps_2
H\exip{|\omega_0|}.
\end{equation}

Then, by Theorem~\ref{t_stat_measure}, we can write
\begin{equation}
\label{mean_numb_part}
 \M(a,b) =
\lambda|\Sph^{d-1}|^{-1}\intl_a^bd\alpha\intl_{\omega_\alpha}du
\intl_{\Sph^{d-1}}dh\,
 \Po^{(\alpha,u),h}[\wp_{-\alpha}(X\cdot\e)
     <\wp_{H-\alpha}(X\cdot\e)].
\end{equation}

Now, let us prove that the rescaled density gradient is given by 
$\vartheta = \lambda\exip{|\omega_0|}$.

\medskip
\noindent
\textit{Proof of Theorem~\ref{t_dens_grad}.} 
Fix an arbitrary $\eps'>0$ and suppose that~$m$ is a (large)
integer. 
Consider the quantity
$t_{m^{-2},m^{-2}}$ defined by~\eqref{prob_good_point}, and suppose
that $H\geq m t_{m^{-2},m^{-2}}^{1/2}$
is large enough to assure that (recall~\eqref{good_points})
\begin{equation}
\label{mnogo_hor}
|\Sph^{d-1}|^{-1} \big|\{(z,h)\in \hD\times \Sph^{d-1} :
   T_\omega^{\eps_1}(z,h)>t_{m^{-2},m^{-2}}\}\big| \leq 2m^{-2}H
     \exip{|\omega_0|}.
\end{equation}

Abbreviate $\phi:=(H/m)^2$ and consider any integer $j\in [1,m]$. 
Suppose that $h\in \Sph^{d-1}$, $z\in\hD$ are such that $z\cdot\e\in
[\frac{(m-j)H}{m},\frac{(m-j+1)H}{m}]$,
and $T_\omega^{m^{-2}}(z,h)\leq t_{m^{-2},m^{-2}}$. Then, since
$\frac{H}{m}\geq t_{m^{-2},m^{-2}}^{1/2}$, from~\eqref{def_Teps} we
obtain that
\begin{align}
 \Po^{z,h}[\wp_{-z\cdot\e}(X\cdot\e)<\wp_{H-z\cdot\e}(X\cdot\e)] 
   &\leq \Po^{z,h}[\wp_{-(m-j)}({\hat Z}^{(\phi)})<\wp_j({\hat
Z}^{(\phi)})] \nonumber\\
   &= \Eo^{z,h} \big(1-R_{m-j,j}({\hat Z}^{(\phi)})\big)\nonumber\\
 &\leq \frac{j}{m} + m^{-2}\nonumber\\
 &\leq \frac{j+1}{m},\label{upper_b_gambler}
\end{align}
and
\begin{align}
 \Po^{z,h}[\wp_{-z\cdot\e}(X\cdot\e)<\wp_{H-z\cdot\e}(X\cdot\e)] 
 &\geq \Po^{z,h}[\wp_{-(m-j+1)}({\hat Z}^{(\phi)})<\wp_{j-1}({\hat
Z}^{(\phi)})] \nonumber\\
   &= \Eo^{z,h} \big(1-R_{m-j+1,j-1}({\hat Z}^{(\phi)})\big)
\nonumber\\
 &\geq \frac{j-1}{m} - m^{-2}\nonumber\\
 &\geq \frac{j-2}{m}. \label{lower_b_gambler}
\end{align}

Also, by the Ergodic Theorem, we can choose~$H$ large enough so that
for all $j=1,\ldots,m$
\begin{equation}
\label{erg_pieces}
 \Big|\frac{m}{H}\big|\hD\cap
  \Xi_{[\frac{(j-1)H}{m},\frac{jH}{m}]
}\big|
  -\exip{|\omega_0|}\Big|
=
\Bigg|\frac{m}{H}\intl_{\frac{(j-1)H}{m}}^{\frac{jH}{m}}
|\omega_\alpha|\,d\alpha 
            - \exip{|\omega_0|}\Bigg| \leq m^{-1}.
\end{equation} 
So, by~\eqref{mean_numb_part}, \eqref{mnogo_hor},
\eqref{lower_b_gambler}, \eqref{erg_pieces}, 
\begin{align}
 \lefteqn{\frac{m}{H}\M\Big(\frac{(m-j)H}{m},
  \frac{(m-j+1)H}{m}\Big)} \nonumber\\
 &\geq
\lambda\frac{m}{H}\times\frac{j-2}{m}\frac{H}{m}
 (\exip{|\omega_0|}-m^{-1}-2m^{-1}\exip{|\omega_0|})\nonumber\\
 &\geq \lambda\frac{j}{m}\exip{|\omega_0|} 
   - \lambda m^{-1}(1+4\exip{|\omega_0|}).
\label{zwei}
\end{align}
Analogously, using~\eqref{upper_b_gambler} instead
of~\eqref{lower_b_gambler}, we obtain
\begin{align}
\lefteqn{
 \frac{m}{H}\M\Big(\frac{(m-j)H}{m},\frac{(m-j+1)H}{m}\Big)}
\nonumber\\
 &\leq \lambda\frac{m}{H}\times 2m^{-2}H\exip{|\omega_0|} + 
\lambda\frac{m}{H}\times\frac{j+1}{m}\frac{H}{m}(\exip{|\omega_0|}
+m^{-1}
)\nonumber\\
 &\leq \lambda\frac{j}{m}\exip{|\omega_0|} 
  + \lambda m^{-1}(2+3\exip{|\omega_0|}).
\label{eins}
\end{align}
Then, we obtain~\eqref{eq_dens_grad} from~\eqref{zwei}
and~\eqref{eins},
and so the proof of Theorem~\ref{t_dens_grad} is concluded.
\qed

At this point, let us formulate an additional result which will be
used in Section~\ref{s_pr_perm}. 
\begin{prop}
\label{p_1/6}
 Define
\begin{align}
\M^*_j &= 
\lambda|\Sph^{d-1}|^{-1}\intl_{\frac{(j-1)m}{H}}^{\frac{jm}{H}}
d\alpha \intl_{\omega_\alpha}du
\intl_{\Sph^{d-1}}dh\,
 \Po^{(\alpha,u),-h}[\wp_{-\alpha}(X\cdot\e)
     <\wp_{H-\alpha}(X\cdot\e)]\nonumber\\
 &\qquad\qquad\qquad\qquad\qquad
  \times \Po^{(\alpha,u),h}[\wp_{H-\alpha}(X\cdot\e)
     <\wp_{-\alpha}(X\cdot\e)],
\label{def_M*}
\end{align}
and suppose that the quenched invariance principle holds.
Then, for any $\eps'>0$ there exists~$m$ such that $\IP$-a.s.
\begin{equation}
\label{eq_1/6}
 \limsup_{H\to\infty}
\max_{j=1,\ldots,m}\Big|
\frac{\M^*_j}{H/m} - \frac{\lambda(j-1/2)(m-j+1/2)}{m}
\exip{|\omega_0|}\Big| 
        < \eps'.
\end{equation}
\end{prop}
\textit{Proof.} The proof is quite analogous to the proof of
Theorem~\ref{t_dens_grad}.
\qed

\medskip

Now, we calculate the limiting rescaled current.

\medskip
\noindent
\textit{Proof of Theorem~\ref{t_current}.}
First, we obtain an upper and a lower bounds for ${\tilde
G}_{i,j,m}$, where $i+j=m$.
By e.g.\ the formula \textbf{1}.2.0.2 of~\cite{BM}, we have
\[
 P[\wp_a(B^{({\hat\sigma})})\leq t] =
\intl_0^t\frac{|a|}{\sqrt{2\pi}{\hat\sigma}s^{3/2}}
   \exp\Big(-\frac{a^2}{2{\hat\sigma}^2s}\Big)\, ds.
\]
So, for $i\leq m^{3/5}$
\begin{align}
 {\tilde G}_{i,j,m} &\geq P[\wp_i(B^{({\hat\sigma})})\leq m]
    - P[\wp_{-j}(B^{({\hat\sigma})}) <
\wp_i(B^{({\hat\sigma})})]\nonumber\\
 & \geq -m^{-2/5} + \intl_0^m\frac{i}{\sqrt{2\pi}{\hat\sigma}s^{3/2}}
   \exp\Big(-\frac{i^2}{2{\hat\sigma}^2s}\Big)\, ds.
\label{Gtil_lower}
\end{align}
Also, for any $i=1,\ldots,m$,
\begin{align}
 {\tilde G}_{i,j,m} &\leq P[\wp_i(B^{({\hat\sigma})})\leq
m]\nonumber\\
  &= \intl_0^m\frac{i}{\sqrt{2\pi}{\hat\sigma}s^{3/2}}
   \exp\Big(-\frac{i^2}{2{\hat\sigma}^2s}\Big)\, ds.
\label{Gtil_upper}
\end{align}
In particular, for $i>m^{3/5}$, we obtain after some elementary
computations 
that there exists a positive constant~$\gamma'$ such that
\begin{equation}
\label{Gtil_hvost}
{\tilde G}_{i,j,m} \leq \gamma'
m^{1/10}\exp\Big(-\frac{m^{1/5}}{2{\hat\sigma}^2}\Big).
\end{equation}

Next, we employ the same strategy as in the proof of
Theorem~\ref{t_dens_grad}.
Fix a large~$m$, and suppose that $H\geq mt_{m^{-2},m^{-2}}^{1/2}$ 
is such that~\eqref{mnogo_hor} holds. 

Now, let~$Y$ be the expected number of particles that were absorbed
in~$\hDl$ up
to time~$H^2/m$, in the stationary regime. 
Clearly, we have then $J^\omega_H = \frac{\Eo Y}{H^2/m}$. So, one can
write 
\begin{align}
 \Eo Y =& \lambda |\Sph^{d-1}|^{-1}\intl_0^H
d\alpha\intl_{\omega_\alpha}du \intl_{\Sph^{d-1}}dh\,
  \Po^{(\alpha,u),-h}[\wp_{-\alpha}(X\cdot\e) > 
\wp_{H-\alpha}(X\cdot\e)]\nonumber\\
   & \qquad {}\times \Po^{(\alpha,u),h}\Big[\wp_{H-\alpha}
(X\cdot\e)\leq \frac{H^2}{m}, 
           \wp_{H-\alpha}(X\cdot\e) < \wp_{-\alpha}
(X\cdot\e)\Big]\nonumber\\
& + \Eo {\widetilde W}_{H,m},\label{calc_current}
\end{align}
where~${\widetilde W}_{H,m}$ is the mean number of particles that were injected 
in~$\hDl$,  successfully crossed the tube, 
and then hit~$\hDr$ before time~$H^2/m$. 

Suppose that $z,h$ are such that $z\cdot\e\in
[\frac{(m-j)H}{m},\frac{(m-j+1)H}{m}]$,
 $T_\omega^{m^{-2}}(z,h)\leq t_{m^{-2},m^{-2}}$,
$T_\omega^{m^{-2}}(z,-h)\leq t_{m^{-2},m^{-2}}$.
Then, analogously to~\eqref{upper_b_gambler}
and~\eqref{lower_b_gambler}, we write
\begin{align}
 \Po^{z,-h}[\wp_{-z\cdot\e}(X\cdot\e)<\wp_{H-z\cdot\e}(X\cdot\e)] &
\geq \frac{j-2}{m},
  \label{lower_gamb_-h}\\
 \Po^{z,-h}[\wp_{-z\cdot\e}(X\cdot\e)<\wp_{H-z\cdot\e}(X\cdot\e)] &
\leq \frac{j+1}{m}.
  \label{upper_gamb_-h}
\end{align}
Moreover, by~\eqref{Gtil_lower}, for $j\leq m^{3/5}$ (recall that
$\phi=(H/m)^2$),
\begin{align}
 \lefteqn{\Po^{z,h}\Big[\wp_{H-z\cdot\e}(X\cdot\e)\leq
\frac{H^2}{m}, \wp_{H-z\cdot\e}(X\cdot\e) < \wp_{-z\cdot
\e}(X\cdot\e)\Big]} \nonumber\\
 &\geq \Po^{z,h}[\wp_j({\hat Z}^{(\phi)}) \leq m, 
\wp_j({\hat Z}^{(\phi)}) < \wp_{-(m-j)}({\hat Z}^{(\phi)})]
\nonumber\\
 &= \Eo^{z,h} G_{j,m-j,m}({\hat Z}^{(\phi)}) \nonumber\\
 &\geq -m^{-2/5} - m^{-2} +
\intl_0^m\frac{j}{\sqrt{2\pi}{\hat\sigma}s^{3/2}}
   \exp\Big(-\frac{j^2}{2{\hat\sigma}^2s}\Big)\, ds.
\label{capture_lower}
\end{align}
Using~\eqref{Gtil_upper}, we obtain
\begin{align}
\lefteqn{\Po^{z,h} \Big[\wp_{H-z\cdot\e}(X\cdot\e)\leq
\frac{H^2}{m}, \wp_{H-z\cdot\e}(X\cdot\e) < \wp_{-z\cdot
\e}(X\cdot\e)\Big] }~~~ \nonumber\\
 &\leq \Eo^{z,h} G_{j-1,m-j+1,m}({\hat Z}^{(\phi)}) \nonumber\\
 &\leq m^{-2} + \intl_0^m\frac{j}{\sqrt{2\pi}{\hat\sigma}s^{3/2}}
   \exp\Big(-\frac{j^2}{2{\hat\sigma}^2s}\Big)\, ds.
\label{capture_upper}
\end{align}
Thus, using~\eqref{mnogo_hor}, \eqref{erg_pieces},
\eqref{calc_current}, \eqref{lower_gamb_-h},
\eqref{capture_upper}, we obtain for some $C_1,C_2,C_3>0$ (observe
that, in comparison to~\eqref{zwei},
to estimate the product of probabilities in~\eqref{calc_current}, we
have to assume that both
$T_\omega^{m^{-2}}(z,h)$ and $T_\omega^{m^{-2}}(z,-h)$ are less than
or equal to $t_{m^{-2},m^{-2}}$)
\begin{align}
 HJ^\omega_H &= \frac{m}{H} \Eo Y \nonumber\\
 &\geq \lambda \frac{m}{H} \sum_{j\leq m^{3/5}}\frac{j-2}{m} \times
\frac{H}{m}(\exip{|\omega_0|}-m^{-1}-4m^{-1}\exip{|\omega_0|})
         \nonumber\\
  & \qquad\qquad {}\times
\Big(-2m^{-2/5}+\intl_0^m\frac{j}{\sqrt{2\pi}{\hat\sigma}s^{3/2}}
   \exp\Big(-\frac{j^2}{2{\hat\sigma}^2s}\Big)\, ds\Big)\nonumber\\
 &\geq \lambda\exip{|\omega_0|} \sum_{j\leq m^{3/5}}\!
\Big(\frac{j}{m}\intl_0^m\frac{j}{\sqrt{2\pi}{\hat\sigma}s^{3/2}}
   \exp\!\Big(-\frac{j^2}{2{\hat\sigma}^2s}\Big) ds -\nonumber\\
& \qquad\qquad\qquad\qquad\qquad\qquad\qquad\qquad
C_1\frac{j}{m}m^{-2/5}-C_2m^{-1}\Big)\nonumber\\
 &\geq -C_3m^{-1/5} + \lambda\exip{|\omega_0|} \sum_{j\leq
m^{3/5}}\frac{j}{m}
 \intl_0^m\frac{j}{\sqrt{2\pi}{\hat\sigma}s^{3/2}}
   \exp\Big(-\frac{j^2}{2{\hat\sigma}^2s}\Big)\, ds.
\label{curr_lower}
\end{align}

To obtain the corresponding upper bound, fix an arbitrary~$\eps>0$
and suppose that~$m$ is large enough so that~\eqref{eq_cross} of
Lemma~\ref{l_cross} holds for those~$\eps,m$. The term~$\Eo {\widetilde W}_{H,m}$
 of~\eqref{calc_current} can be estimated in the following way:
\[
 \Eo {\widetilde W}_{H,m} \leq C_4 \frac{H^2}{m} \Po^{\hDr}[\CC_H, \TT_H\leq 
        m^{-1}H^2] \leq C_4 \frac{H^2}{m}\times\frac{\eps}{H}, 
\]
so $\frac{m}{H}\Eo {\widetilde W}_{H,m}\leq C_4\eps$. Then,
analogously to~\eqref{curr_lower}, using also~\eqref{Gtil_hvost}, we
have for some $C_5,C_6>0$
\begin{align}
 HJ^\omega_H &\leq \lambda\frac{m}{H}
  \times 4m^{-2}H\exip{|\omega_0|} + C_4\eps
\nonumber\\ 
 & \quad + \lambda\frac{m}{H}
  \sum_{j\leq m^{3/5}} \frac{j+1}{m} \nonumber\\ 
 & \quad\qquad
 \times\frac{H}{m}(\exip{|\omega_0|}+m^{-1})
\Big(m^{-2} + \intl_0^m\frac{j}{\sqrt{2\pi}{\hat\sigma}s^{3/2}}
   \exp\Big(-\frac{j^2}{2{\hat\sigma}^2s}\Big)\, ds \Big) 
  \nonumber\\
 & \quad + \lambda\frac{m}{H} \times
  (m-m^{3/5}) \nonumber\\
& \quad\qquad\times \frac{H}{m}(\exip{|\omega_0|}+m^{-1})
\Big(m^{-2} 
        +
\gamma'm^{1/10}\exp\Big(-\frac{m^{1/5}}{2{\hat\sigma}^2}
\Big)\Big)\nonumber\\
 &\leq \lambda\exip{|\omega_0|} \sum_{j\leq
m^{3/5}}\frac{j}{m}\intl_0^m\frac{j}{\sqrt{2\pi}{\hat\sigma}s^{3/2}}
   \exp\Big(-\frac{j^2}{2{\hat\sigma}^2s}\Big)\, ds \nonumber\\
 & \quad + C_4\eps + C_5m^{-1} +
C_6m^{1/10}\exp\Big(-\frac{m^{1/5}}{2{\hat\sigma}^2}\Big).
\label{curr_upper}
\end{align}

Now, observe that
\begin{align*}
\lefteqn{\lim_{m\to\infty} \sum_{j\leq
m^{3/5}}\frac{j}{m}\intl_0^m\frac{j}{\sqrt{2\pi}{\hat\sigma}s^{3/2}}
   \exp\Big(-\frac{j^2}{2{\hat\sigma}^2s}\Big)\, ds}\\
 &= \lim_{m\to\infty} \sum_{j\leq m^{3/5}}\frac{j}{m}
     \intl_0^1\frac{j}{\sqrt{2\pi}{\hat\sigma}m^{3/2}s^{3/2}}
   \exp\Big(-\frac{j^2}{2{\hat\sigma}^2ms}\Big)m\, ds\\
 &= \lim_{m\to\infty} \sum_{j\leq m^{3/5}}
 \frac{1}{\sqrt{m}}\intl_0^1\frac{(j/\sqrt{m})^2}{\sqrt{2\pi}{
\hat\sigma}s^{3/2}}
   \exp\Big(-\frac{(j/\sqrt{m})^2}{2{\hat\sigma}^2s}\Big)\, ds\\
 &= \intl_0^{\infty}dr \intl_0^1 ds\,
\frac{r^2}{\sqrt{2\pi}{\hat\sigma}s^{3/2}}
   \exp\Big(-\frac{r^2}{2{\hat\sigma}^2s}\Big) \\
 &= \intl_0^1 \frac{ds}{s}\intl_0^{\infty}dr\,
\frac{r^2}{\sqrt{2\pi}{\hat\sigma}s^{1/2}}
   \exp\Big(-\frac{r^2}{2{\hat\sigma}^2s}\Big)\\
 &= \intl_0^1 \frac{ds}{s} \times \frac{{\hat\sigma}^2s}{2}\\
 &= \frac{{\hat\sigma}^2}{2}.
\end{align*}
With this observation, Theorem~\ref{t_current} 
follows from~\eqref{curr_lower} and~\eqref{curr_upper}.
\qed

\subsection{Proof of Theorems~\ref{t_exp_perm} and~\ref{t_perm}}
\label{s_pr_perm}
Observe that, since the particles are independent, 
the Knudsen gas in the finite tube~$\tD$ can be
regarded as a $M/G/\infty$ queueing system; moreover, using e.g.\
Theorem~2.1 of~\cite{CPSV1} it is straightforward to obtain that the
service time (which is the lifetime of a newly injected particle) is
a random variable with exponential tail. Then, let us recall the
following basic identity of queuing theory (known as Little's
theorem): 
\begin{prop}
\label{prop_Little}
Suppose that $\Lambda_a$ is
the arrival rate, $q$ is the mean number of customers in the system,
and~$T$ is the mean time a customer spends in the system, then
$T=q/\Lambda_a$.
\end{prop}

\noindent
\textit{Proof.} See e.g.\ Section~5.2 of~\cite{C81}. To understand
intuitively why this fact holds true,
 one may reason in the following way: 
by large time~$t$, the total time of all the
customers in the system would be (approximately)~$qt$ on one hand,
and~$T\Lambda_a t$ on the other hand.
\qed

\medskip
\noindent
\textit{Proof of Theorem~\ref{t_exp_perm}.}
This result almost immediately follows from
Theorem~\ref{t_dens_grad} by using Proposition~\ref{prop_Little}.
First, for the gas of independent particles the arrival rate is
\begin{equation}
\label{def_Lambda_a}
 \Lambda_a =
\frac{\lambda|{\tilde\omega}_0|}{\gamma_d|\Sph^{d-1}|},
\end{equation}
recall that the particles are injected in~$\tDl$ only.
Then, from Theorem~\ref{t_dens_grad} it is straightforward to obtain
that for the mean number of particles~$q_H$ in the system, we have
\[
 \lim_{H\to\infty} \frac{q_H}{H} = 
      \frac{\lambda H\exip{|\omega_0|}}{2}.
\]
Then, Proposition~\ref{prop_Little} implies~\eqref{q_exp_perm}. To
prove the corresponding annealed result, note that $q_H\leq \lambda
H|\Xi|$ by
Theorem~\ref{t_stat_measure}~(ii). So, applying the bounded
convergence
theorem, we obtain~\eqref{a_exp_perm}.
\qed

\medskip
\noindent
\textit{Proof of Theorem~\ref{t_perm}.}
First, observe that in the stationary regime the particles
leave the system at the right boundary with rate~$J_H$, and this
should be equal
to the entrance rate~$\Lambda_a\Po[\CC_H]$ of the particles which
cross the tube, with~$\Lambda_a$ from~\eqref{def_Lambda_a}.
 So, \eqref{q_perm_prob} follows from Theorem~\ref{t_current}.

To prove~\eqref{a_perm_prob}, observe that, by
using Lemma~\ref{l_escape_upper} with~$B=\tDl$
and~$F=\hDr$, we obtain that for some positive
constants~$C_1,C_2$ which do not depend on~$\omega$
\[
 H\Po[\CC_H] \leq C_1 + \frac{C_2}{H}
   \intl_{{\tilde F}^\omega(0,H)}b(x)\, d\nuo(x).
\]
By~\eqref{conv_b}, the collection of random variables
$(H\Po[\CC_H],H>1)$ is uniformly integrable, and this
implies~\eqref{a_perm_prob}.

In order to prove~\eqref{q_cond_crossing}, denote by~$q'_H$
the mean number of particles in the stationary regime that
\emph{will exit} at~$\hDr$. Observe that, by 
Theorem~\ref{t_stat_measure}~(ii) and Proposition~\ref{p_1/6},
\[
 \lim_{H\to\infty} \frac{q'_H}{H} = \lambda\exip{|\omega_0|}
       \intl_0^1 x(1-x)\, dx = \frac{\lambda\exip{|\omega_0|}}{6}.
\]
So, using~\eqref{q_perm_prob} and Proposition~\ref{prop_Little},
we obtain~\eqref{q_cond_crossing}. The relations~\eqref{q_crossing}
and~\eqref{q_not_crossing} follow from~\eqref{q_cond_crossing}
and~\eqref{q_perm_prob}.

Now, observe that~\eqref{a_crossing} and~\eqref{a_not_crossing}
immediately follow from~\eqref{q_crossing}, \eqref{q_not_crossing},
and~\eqref{a_exp_perm}, so now it remains only 
to prove~\eqref{a_cond_crossing}. Let~$\sigma_1:=\tau^+(\tDl)$,
$\sigma_{k+1}=\min\{m>\sigma_k : \xi_m\in\tDl\}$ be the moments
of successive visits to~$\tDl$ for the process in the finite tube.
By Corollary~\ref{c_trans_dens}, $\xi_{\sigma_k}$ is uniformly
distributed in~$\tDl$ for all~$k$, and so we can write
\begin{equation}
\label{during_k}
 \Po^{\tDl}[\tau(\hDr)<\sigma_k] \leq k\Po[\CC_H].
\end{equation}
Then, using~\eqref{during_k}, Lemma~\ref{l_Dirichlet_lower},
 and the fact that the random variables
$(Z_j,j\geq 1)$ are independent of everything, we obtain
\begin{align*}
 \frac{C_3}{H} &\leq \Po[{\hat\CC}_H]\\
  &\leq \sum_{k=1}^\infty 
    \Po[{\hat\tau}(\hDr)<{\hat\tau}(\tDl),
       \sigma_{k-1}<Z_1+\cdots+Z_{{\hat\tau}(\hDr)}<\sigma_k]\\
 & \leq \Po[Z_1+\cdots+Z_j\neq \sigma_\ell
 \text{ for all $\ell<k$ and all~$j$ }\mid \tau(\hDr)<\sigma_k]\\
 & \qquad \times \Po[\tau(\hDr)<\sigma_k]\\
 & \leq \Po[\CC_H] \sum_{k=1}^\infty 
   k(1-N^{-1})^{\lceil\frac{k-1}{N}\rceil},
\end{align*}
and this implies that $\Po[\CC_H]\geq C_4/H$ for some $C_4>0$
not depending on~$\omega$. Since $q'_H\leq \lambda H|\Xi|$,
one obtains~\eqref{a_cond_crossing} from the bounded convergence
theorem.
\qed

\section*{Acknowledgements}
We thank Takashi Kumagai for pointing us reference~\cite{C}.
The work
of F.C.\ was partially supported by CNRS (UMR 7599
``Probabilit{\'e}s et Mod{\`e}les Al{\'e}atoires'') and ANR
Polintbio. 
S.P.\
was partially supported by  CNPq (300886/2008--0). G.M.S.\ thanks
DFG (Priority programme SPP 1155) for financial support.
The work of M.V.\ was partially supported 
         by CNPq (304561/2006--1).
S.P.\ and M.V.\ also thank FAPESP (2009/52379--8), CNPq
(471925/2006--3, 472431/2009--9),
 and CAPES/DAAD (Probral) for financial support.

\end{document}